\documentclass[11pt]{amsart} \hoffset -1cm \voffset -1in
\usepackage{marvosym}
\usepackage{amsmath}
\usepackage{amssymb}
\usepackage{amsthm}
\usepackage{color}
\usepackage{graphicx}
\usepackage{array}
\usepackage{float}
\usepackage{phonetic}
\usepackage{pst-all}
\usepackage{pb-diagram}
\usepackage[T1]{fontenc}
\usepackage{ stmaryrd }
\usepackage{url}

\newtheorem{prop}{Proposition}[section]
\newtheorem{lema}[prop]{Lemma}
\newtheorem{teor}[prop]{Theorem}
\newtheorem{hecho}[prop]{Fact}
\newtheorem{obs}[prop]{Observation}
\newtheorem{preg}[prop]{Question}
\newtheorem{corol}[prop]{Corollary}
\newtheorem{conj}[prop]{Conjecture}

\theoremstyle{definition}
\newtheorem{definicion}[prop]{Definition}
\newtheorem{ejem}[prop]{Example}
\newtheorem{rmk}[prop]{Remark}
\newtheorem{const}[prop]{Construction}

\def\Ind#1#2{#1\setbox0=\hbox{$#1x$}\kern\wd0\hbox to 0pt{\hss$#1\mid$\hss}     
\lower.9\ht0\hbox to 0pt{\hss$#1\smile$\hss}\kern\wd0}
\def\Notind#1#2{#1\setbox0=\hbox{$#1x$}\kern\wd0\hbox to 0pt{\mathchardef
\nn=12854\hss$#1\nn$\kern1.4\wd0\hss}\hbox to
0pt{\hss$#1\mid$\hss}\lower.9\ht0 \hbox to
0pt{\hss$#1\smile$\hss}\kern\wd0}

\def\ind{\mathop{\mathpalette\Ind{}}}                       
\def\thind{\mathop{\mathpalette\Ind{}}^{\text{\th}} }       
\def\uth{\text{U}^{\text{\th}} }                            
\newcommand{\bp}{\begin{prop}}
\newcommand{\ep}{\end{prop}}
\newcommand{\bd}{\begin{definicion}}
\newcommand{\ed}{\end{definicion}}
\newcommand{\bej}{\begin{ejem}}
\newcommand{\eej}{\end{ejem}}
\newcommand{\bl}{\begin{lema}}
\newcommand{\el}{\end{lema}}
\newcommand{\bh}{\begin{hecho}}
\newcommand{\eh}{\end{hecho}}
\newcommand{\bpreg}{\begin{preg}}
\newcommand{\epreg}{\end{preg}}
\newcommand{\bo}{\begin{obs}}
\newcommand{\eo}{\end{obs}}
\newcommand{\bcon}{\begin{conj}}
\newcommand{\econ}{\end{conj}}
\newcommand{\brmk}{\begin{rmk}}
\newcommand{\ermk}{\end{rmk}}
\newcommand{\bc}{\begin{corol}}
\newcommand{\ec}{\end{corol}}
\newcommand{\bconst}{\begin{const}}
\newcommand{\econst}{\end{const}}
\newcommand{\bdem}{\begin{proof}}
\newcommand{\edem}{\end{proof}}
\newcommand{\benum}{\begin{enumerate}}
\newcommand{\eenum}{\end{enumerate}}
\newcommand{\bitem}{\begin{itemize}}
\newcommand{\eitem}{\end{itemize}}
\newcommand{\be}{\begin{ejem}}
\newcommand{\ee}{\end{ejem}}
\newcommand{\bt}{\begin{teor}}
\newcommand{\et}{\end{teor}}
\newcommand{\Los}{\L o\'s}
\newcommand{\pte}[1]{\textcolor{red}{\textbf{#1}}}

\newcommand{\modulo}[1]{\,\,\left(\operatorname{mod}\,\,#1\right)}

\newcommand{\R}{\mathbb{R}}

\newcommand{\Z}{\mathbb{Z}}

\newcommand{\C}{\mathcal{C}}
\newcommand{\Nat}{\mathbb{N}}
\newcommand{\Q}{\mathbb{Q}}

\newcommand{\ov}[1]{\overline{#1}}

\newcommand{\ci}[1]{\text{\textcircled{#1}}}
\title{Ordered asymptotic classes of finite structures}
\author{Dar\'io Garc\'ia\\ \\ University of Leeds}
\address{Dar\'io Garc\'ia  \\ School of Mathematics \\ University of Leeds \\ Leeds, UK.}
\curraddr{Departamento de Matemáticas. Universidad de los Andes\\ Carrera 1 No. 18A-10, Edificio H, Bogot\'a 111711, Colombia.}
\email{da.garcia268@uniandes.edu.co}

\begin{document}
\begin{abstract} We introduce the concept of o-asymptotic classes of finite structures, melding ideas coming from 1-dimensional asymptotic classes and o-minimality. Along with several examples and non-examples of these classes, we present some classification theory results of their infinite ultraproducts: Every infinite ultraproduct of structures in an o-asymptotic class is superrosy of $\uth$-rank 1, and NTP$_2$ (in fact, inp-minimal).
\end{abstract}
\maketitle
\setlength{\parindent}{0pt}

\section{Introduction}
In \cite{MS}, Macpherson and Steinhorn develop the notion of 1-dimensional asymptotic classes, which are classes of finite structures with notions of measure and dimension coming from the study of the size of their definable sets. Specifically, they define

\bd Let $L$ be a first order language, and $\mathcal{C}$ be a collection of finite $L$-structures. Then $\mathcal{C}$ is a \emph{1-dimensional asymptotic class} if for every formula $\varphi(x,\ov{y})$, with $\ov{y}=(y_1,\cdots,y_m)$, the following conditions hold:
\benum
\item There is a positive constant $C$ and a finite set $E\subseteq \R^{>0}$ such that for every $M\in\mathcal{C}$ and $\ov{a}\in M^m$, either $|\varphi(M,\ov{a})|\leq C$, or for some $\mu\in E$,
\[\left||\varphi(M,\ov{a})|-\mu|M|\right|\leq C|M|^{1/2}.\]
\item For every $\mu\in E$, there is an $L$-formula $\varphi_\mu(\ov{y})$ such that for all $M\in\C$, $\varphi_\mu(M^m)$ is precisely the set of tuples $\ov{a}\in M^m$ with
\[\left||\varphi(M,\ov{a})|-\mu|M|\right|\leq C|M|^{1/2}.\]
\eenum
\ed

The seminal example of these classes is the class of finite fields, for which the conditions above appear as a remarkable theorem of Chatzidakis, van den Dries and Macintyre  (see \cite{ChVdDMac}).  With this definition, which is a condition on the definable sets in only one variable, Macpherson and Steinhorn obtain results about the control of the size of definable sets in several variables as well as results concerning the behavior of infinite ultraproducts of structures in such classes. For instance, they show that if every ultraproduct in a class $\C$ is strongly minimal, then $\C$ is a 1-dimensional asymptotic class. Furthermore, they prove that every ultraproduct of a 1-dimensional class is supersimple of $U$-rank 1.

An easy example of a class of finite structures which is not a 1-dimensional class is the class of all finite totally ordered sets, which fails property (1) because the formula $x<a$ can pick out an arbitrary proper initial segment of a structure as the parameter $a$ varies. However, the only definable sets in one variable on the structures of this class (and in their ultraproducts) are finite unions of intervals and points, implying that the structures involved are \emph{o-minimal}.

O-minimality and its variants are conditions for theories of infinite ordered structures that provide good structural properties of their models. Our aim here is to isolate conditions on classes of finite linearly ordered structures to get nice asymptotic properties, melding ideas of asymptotic classes and o-minimality.

With this idea in mind, we propose a definition of \emph{o-asymptotic classes} as an adaptation of the definition of 1-dimensional asymptotic classes in the context of totally ordered structures.\\

This paper is organized as follows: in Section 2 we present the definition of o-asymptotic classes and some basic properties and examples. In Section 3 we start the study of infinite ultraproducts of structures in o-asymptotic classes, for which our main interest is to place them in the map of classification theory. Among the main results, we have that if every ultraproduct of a class of finite totally ordered structures is o-minimal, then the class is o-asymptotic (Proposition 3.1), and furthermore, every ultraproduct of an o-asymptotic class is NTP$_2$ and super rosy of $\uth$-rank 1. (Theorems 3.5 and 3.6)

In Section 4 we present our main example of an o-asymptotic class, that consists of the class of cyclic groups $(\Z/(2N+1)\Z,+)$ with a natural order on the equivalence classes given by $-\ov{N}<\cdots <-\ov{1}<\ov{0}<\ov{1}<\cdots<\ov{N}$. Section 5 is devoted to presenting several non-examples of o-asymptotic classes, including several classes of ordered graphs. Finally, some quantifier elimination results needed for the examples are included in the Appendix, as well as a proof of quasi-o-minimality for infinite ultraproducts of certain class of finite linearly ordered structures.

\subsection{Rosy theories and theories with NTP$_2$}

During the last few years, the class of NTP$_2$ theories (which contains both simple theories and theories with NIP) has received the attention of model-theorists, especially after the results in \cite{Che,CheKaSi,CheHi}. Even though NTP$_2$ it was defined in \cite{Sh1,Sh80}, its systematic study was initiated in \cite{Che}. Now we give the definitions of these theories.

\bd \ \benum
\item An \emph{inp-pattern of depth $\kappa$} in a model $M$ consists on a sequence of formulas $\langle \varphi_\alpha(\ov{x},\ov{y}_{\alpha}):\alpha<\kappa\rangle$, integers $k_\alpha$ and tuples $\langle \ov{a}_{\alpha,i}:\alpha<\kappa, i<\omega\rangle$ from $M$ such that:
\bitem
\item For each $\alpha<\kappa$, the set $\{\varphi_\alpha(\ov{x},\ov{a}_{\alpha,i}):i<\omega\}$ is $k_\alpha$-inconsistent.
\item For any function $f:\kappa\to\omega$, the set $\{\varphi_\alpha(\ov{x},\ov{a}_{\alpha,f(\alpha)}):\alpha<\kappa\}$ is consistent.
\eitem
\item We say that a formula $\varphi(\ov{x},\ov{y})$ has \emph{TP$_2$}  (or has the \emph{tree property of the second kind}) relative to $M$ if there is an inp-pattern of depth $\omega$ in $M$ with $\varphi_\alpha(\ov{x},\ov{y}_\alpha)=\varphi(\ov{x},\ov{y})$.
\item We say that $M$ \emph{is NTP$_2$} if there is no $L$-formula witnessing TP$_2$ in any model of $\operatorname{Th}(M)$. A theory $T$ has NTP$_2$ if every model of $T$ is NTP$_2$.
\item We say that a structure $M$ is \emph{inp-minimal} if there is no inp-pattern of depth $2$ with formulas in one object variable (i.e., $|\ov{x}|=1$) in any model of $\operatorname{Th}(M)$.
\eenum
\ed

A way to understand inp-patterns is to think about an ``infinite rectangular array'' of formulas $\varphi_\alpha(\ov{x},\ov{a}_{\alpha,i})$ (where the index $i$ varies on the rows while $\alpha$ varies on the columns) for which every row is $k_\alpha$-inconsistent, but every ``descending path'' is consistent. \\

There are recent results providing examples of natural theories which have NTP$_2$ but are neither simple theories nor theories with NIP: they include ultraproducts of $p$-adics (see \cite{CheKaSi}) and the theory of the non-standard Frobenius automorphism acting on an algebraically closed valued field of equicharacteristic $0$ (see \cite{CheHi}). 

\bh If $T$ has TP$_2$, there is some formula $\varphi(x,\ov{y})$ with $|x|=1$ that has TP$_2$. Thus, in particular, if a structure is inp-minimal then it is NTP$_2$.
\eh

\bh [Chernikov, Lemma 2.2 in \cite{Che}] \label{mutuallyindisc}  If there is an array $\langle \varphi_\alpha(x,\ov{b}_{\alpha,i}):\alpha<\kappa, i<\omega\rangle$ witnessing an inp-pattern of depth $\kappa$, then the sequences $I_\alpha=\langle b_{\alpha,i}:i<\omega\rangle$ may be assumed to be mutually indiscernible sequences.
\eh

\subsection{\th-Forking and Rosy Theories}
The notion of \th-forking (read as \emph{thorn-forking}) was introduced in Onshuus's Ph.D. Thesis and appeared as a generalization of forking to contexts in which non-forking failed to provide a nice independence relation, such as that of o-minimal theories.

\bd \benum
\item A formula $\varphi(\ov{x},\ov{b})$ \emph{strongly divides} over $A$ if $\ov{b}$ is not algebraic over $A$ and there is $k<\omega$ such that the set $\{\varphi(\ov{x},\ov{b}'):\ov{b}'\models \operatorname{tp}(\ov{b}/A)\}$ is $k$-inconsistent.
\item We say that $\varphi(\ov{x},\ov{b})$ \emph{\th-divides} over $A$ if there is a finite tuple $e$ (possibly in $M^{eq}$) such that $\varphi(\ov{x},\ov{b})$ strongly divides over $Ae$.
\item We say that $\varphi(\ov{x},\ov{b})$ \emph{\th-forks} over $A$ if there are formulas $\psi_1(\ov{x},\ov{c}_1),\ldots,\psi_n(\ov{x},\ov{c}_n)$ such that each $\psi_i(\ov{x},\ov{c}_i)$ \th-divides over $A$ and $\displaystyle{\varphi(\ov{x},\ov{b})\vdash \bigvee_{i=1}^n \psi_i(\ov{x},\ov{c}_i)}.$
\item A type \emph{\th-forks (\th-divides) over $A$} if it implies a formula that \th-forks (\th-divides) over $A$.
\item We write $\ov{a}\thind_A B$ (read as \emph{$\ov{a}$ is thorn-independent of $B$ over $A$}) to denote that $\operatorname{tp}(\ov{a}/B)$ does not \th-fork over $A$.
\eenum
\ed

Just as simple theories are theories in which forking-independence has nice properties and can be characterized as those theories where forking-independence is symmetric (meaning that $\ov{a}\ind_A \ov{b}$ if and only if $\ov{b}\ind_A\ov{a}$), there is a class of theories called \emph{rosy} in which \th-forking has desirable properties. Rosy theories can can also be characterized as those theories where \th-independence is a symmetric independence relation. The main properties of \th-forking and rosy theories were investigated in \cite{On1}, where Onshuus also showed that all simple theories and all o-minimal theories are rosy.

\bd \label{superrosy} \ 
\benum
\item We say that $q\in S(B)$ is a \emph{\th-forking extension} of $p\in S(A)$ (with $A\subseteq B$) if $q$ is an extension of $p$ and the type $q$ \th-forks over $A$. Otherwise, we called it a non-\th-forking extension of $p$.
\item We define the \emph{$\uth$-rank} (read as U-thorn-rank) to be the foundation rank for \th-forking. Namely, $\uth(p(\ov{x}))\geq 0$ if and only if $p(\ov{x})$ is consistent, $\uth(p(\ov{x}))\geq \alpha+1$ if and only if there is a \th-forking extension $q(\ov{x})$ of $p(\ov{x})$ such that $U(q(\ov{x}))\geq \alpha$, and for a limit ordinal $\lambda$, $\uth(p(\ov{x}))\geq \lambda$ if and only if $\uth(p(\ov{x}))\geq \alpha$ for every $\alpha<\lambda$.
\item A structure $M$ is said to be \emph{superrosy of $\uth$-rank $n$} if there is a $1$-type $p(x)$ such that $\uth(p(x))=n$, but there is no $1$-type $q(x)$ with $\uth(q(x))\geq n+1$, or equivalently, if the maximal length of a \th-dividing chain for $1$-types in $M$ is $n$.
\eenum
\ed

 In all supersimple theories, U-rank and $\uth$-rank coincide. On the other hand, it is known that in the presence of a definable order forking is very different from \th-forking. For example, in the theory $\operatorname{Th}(\Q,<)$ we have that the formula $\varphi(x):=a<x<b$ divides over the empty set (despite the fact that the o-minimal dimension of $\varphi(M)$ is $1$), but it does not \th-fork over the empty set. This behavior generalizes to all o-minimal theories, where \th-independence coincides with the usual notion of independence, and $\uth$-rank corresponds to the o-minimal dimension on definable sets. (cf. \cite[Section 5.2]{On1})

\subsection{A lemma about measure theory} The following fact is a consequence of the truncated inclusion-exclusion principle, and will be used in Section \ref{ultraoclasses}. A proof of it can be found in \cite[Proposition 2.2.10]{GMS}.

\bh \label{kintersections} Let $X$ be a measure space with $\mu(X)=1$ and fix $0<\epsilon\leq \dfrac{1}{2}$. Let $\langle A_i:i<\omega\rangle$ be a sequence of measurable subsets of $X$ such that $\mu(A_i)\geq \epsilon$ for every $i$. Then, for every $k<\omega$ there are $i_1<i_2<\ldots<i_k$ such that \[\mu\left(\bigcap_{j=1}^k A_{i_j}\right)\geq \epsilon^{3^{k-1}}.\]
\eh

\section{o-asymptotic classes and cell decomposition results} \label{odefcelldecomp}
\bd \label{def-weak-o-class} Let $\mathcal{C}$ be a class of finite linearly ordered structures in a language $L$ containing $<$. We say $\mathcal{C}$ is a \emph{\textbf{weak}-o-asymptotic class} if for every formula $\varphi(x;y_1,\ldots,y_m)$ there is a constant $C=C_\varphi>0$ and $k=k_\varphi\geq 1$ and a finite set $E_\varphi\subseteq [0,1]^k$ such that:

\benum
\item For every $M\in\mathcal{C}$ and $\ov{a}\in M^m$ there are elements
\[c_0=\min M\leq c_1\leq \ldots \leq c_k=\max M\] and a tuple $\ov{\mu}\in E$ such that:\\
(*) For  every $i=1,2,\ldots,k$, either
\[\text{\hspace{1cm}}\begin{cases}\mu_i=0 \text{\ \  and\ \ } |\varphi(M,\ov{a})\cap [c_{i-1},c_i]|\leq C\\ \text{\hspace{2cm}or} \\ \mu_i>0 \text{\ \ and\ \ } ||\varphi(M,\ov{a})\cap [c_{i-1},c_i]|-\mu_i|[c_{i-1},c_i]||\leq C|[c_{i-1},c_i]|^{1/2}.\end{cases}\]\item For every $\ov{\mu}\in E$ there is a formula $\varphi_{\ov{\mu}}(\ov{y};z_0,z_1,\ldots,z_{k})$ such that for every $M\in\mathcal{C}$ and $\ov{a}\in M$,
\[M\models \varphi_{\ov{\mu}}(\ov{a};c_1,\ldots,c_{k}) \text{\hspace{0.5cm}implies \hspace{0.5cm} (*) holds}.\]
\eenum
\ed

\bd \label{defoclasses} Let $\mathcal{C}$ be a class of finite linearly ordered structures in a language $L$ containing $<$. We say $\mathcal{C}$ is an \emph{o-asymptotic class} if for every $m\in\Nat$ and formula $\varphi(x;y_1,\ldots,y_m)$ there is a constant $C>0$ and $k\geq 1$ and a finite set $E\subseteq [0,1]^k$ such that:

\benum
\item For every $M\in\mathcal{C}$ and $\ov{a}\in M^m$ there are elements
\[c_0=\min M<c_1<\ldots<c_k=\max M\] and a tuple $\ov{\mu}\in E$ such that:\\
(*) For  every $i=1,2,\ldots,k$, either
\[\text{\hspace{1cm}}\begin{cases}\mu_i=0 \text{\ \  and\ \ } |\varphi(M,\ov{a})\cap [c_{i-1},c_i]|\leq C\\ \text{\hspace{2cm}or} \\ \mu_i>0 \text{\ \textbf{ and for every }$(u,v)\subseteq [c_{i-1},c_i]$,}\\
 \hspace{1cm} ||\varphi(M,\ov{a})\cap [u,v]|-\mu_i|[u,v]||\leq C|[u,v]|^{1/2}.\end{cases}\]
 
 \item (Definability condition) For every $\ov{\mu}\in E$ there is a formula $\varphi_{\ov{\mu}}(\ov{y};z_1,\ldots,z_{k})$ such that for every $M\in\mathcal{C}$,
\[M\models \varphi_{\ov{\mu}}(\ov{a};c_1,\ldots,c_{k}) \text{\hspace{0.5cm}implies \hspace{0.5cm} (*) holds}.\]
\eenum
\ed

\brmk Roughly speaking, a class of finite ordered structures is weakly o-asymptotic if every formula in one variable admits a decomposition into a fixed number of intervals such that on each interval it behaves as in 1-dimensional asymptotic classes.  Being o-asymptotic requires also this decomposition to be \emph{uniform} in the sense that the definable set is uniformly distributed along each of the intervals $[c_{i-1},c_i]$. We believe that these two notions are equivalent, but we have not been able to prove this equivalence. The main difficulty here has been the non-additive nature of the ``error term'' $[u,v]|^{1/2}$. \\

\noindent Throughout this paper, we will focus mainly on the more robust definition of o-asymptotic classes.
\ermk

\brmk
The error term $C|[u,v]|^{1/2}$ has been chosen in analogy to the error term $C |M|^{1/2}$ used for 1-dimensional asymptotic classes (see \cite[Definition 1.2]{MS}), which itself is based on \cite{ChVdDMac} and ultimately on Lang-Weil estimates for algebraic varieties defined over finite fields. For asymptotic classes it is usual to relax the error term, for example requiring it to be $C|M|^{1-\epsilon}$, or even $o(|M|)$ as \cite[Definition 2.1]{Elw}. With this new error term, new examples of asymptotic classes can be found, while all the structural results (classification of ultraproducts, uniform quantifier elimination) are preserved.\\

The same variations can be considered for o-asymptotic classes of finite structures, and the proofs of the main results (Proposition \ref{O-classes-o-minimal}, Theorems \ref{Uthrank1}, \ref{O-NTP2}) can be easily adapted to this more general context. In addition, all the non-examples we will present in Section 5 can also be modified to provide non-examples even with this more general definition. We will however maintain the error term $C|[u,v]|^{1/2}$ to provide explicit calculations.
\ermk

\noindent \textbf{Notation:} Assume we are working in an o-asymptotic class $\mathcal{C}$, and let $M$ be a structure in $\C$. \bitem \item We say that $\varphi(x;\ov{a})$ \emph{admits a decomposition with proportion $\ov{\mu}$} to mean that there are $c_1,\ldots,c_k$ such that condition (1) of Definition \ref{defoclasses} holds for $\varphi(x;\ov{a})$ and $c_1<\cdots <c_k$.
\item When we say \emph{by uniformity of $\mu_i$ in $[c,d]$} or \emph{by the uniformity of distribution of $\varphi(x;\ov{a})$ in $[c,d]$} we mean that, for every $[u,v]\subseteq [c,d]$,
\[||\varphi(M,\ov{a})\cap [u,v]|-\mu_i|[u,v]||\leq C|[u,v]|^{1/2}.\] 
\eitem

\brmk \label{m-functions} Note that the definability condition provides also a definable way to obtain decompositions for definable sets in one variable, uniformly on the finite structures of an o-asymptotic class. Namely, given a formula $\varphi(x,\ov{y})$ we can define the functions $m_{\varphi,\ov{\mu}}^{i}(\ov{y})$ for $i=0,\ldots,k_\varphi$ as follows:
\bitem
\item For every $\ov{y}$, $m_\varphi^0(\ov{y})=\min M$, $m_{\varphi}^{k_\varphi}(\ov{y})=\max M$
\item For $i=1,\ldots,k_\varphi-1$, we define 
\begin{align*}
&M\models m_\varphi^i(\ov{a})=d\\
&\Leftrightarrow d=\min \left\{d'\in M: \exists z_{i+1},\ldots,z_k\,\varphi_{\ov{\mu}}\left(\ov{a};m_{\varphi,\ov{\mu}}^0(\ov{a}),\ldots,m_{\varphi,\ov{\mu}}^{i-1}(\ov{a}),d',z_{i+1},\ldots,z_k\right)\right\}.
\end{align*}
\eitem
The functions $m_{\varphi,\ov{\mu}}^{i}(\ov{y})$ are $\emptyset$-definable, and when evaluated in some parameter $\ov{a}$ they provide the end points of the intervals for a good decomposition of the definable set $\varphi(x,\ov{a})$. These functions will play an important role in Section \ref{ultraoclasses} where we study model-theoretic properties of infinite ultraproducts of o-asymptotic classes.
\ermk

\brmk Suppose $\mathcal{C}$ is a class of discrete linearly ordered $L$-structures with endpoints, with $L$ containing the symbol $\{<\}$. If we define the expanded language $L'=\{<,S,S^{-1},\min,\max\}$, then every structure in $\mathcal{C}$ has a canonical expansion to an $L'$-structure, where the constant symbols $\min,\max$ are interpreted as the minimum and maximum elements (respectively), and the unary functions $S,S^{-1}$ are interpreted as the successor and the predecessor functions. For the sake of completeness, we also put $S(\max)=\max$ and $S^{-1}(\min)=\min$.\\

These canonical expansions will be used several times to obtain uniform quantifier elimination results, which can be seen as a first step towards proving that a class $\mathcal{C}$ of finite linearly ordered structures is an o-asymptotic class.
\ermk

\bej \label{example-Cord} \textbf{The class $\mathcal{C}_{ord}$ of finite linear orders is an o-asymptotic class.}

It is easy to show that this class admits uniform quantifier elimination in the language $L'=\{<,S,S^{-1},\min,\max\}$ (see Lemma \ref{UQEcord}). Moreover, every formula $\varphi(x;\ov{y})$ can be written (uniformly across the class $\mathcal{C}_{ord}$) as a fixed finite disjunction of formulas of the form \[\bigwedge_i \tau_1^i(\ov{y})<x \wedge \bigwedge_j x<\tau^j_2(\ov{y})\wedge \bigwedge_\ell x=\tau^\ell_3(\ov{y}),\] where $\tau_1^i,\tau_2^j,\tau_3^\ell$ are terms depending on the tuple $\ov{y}$ in the language $L'$. Thus, for every tuple $\ov{a}$ from a structure in $\mathcal{C}_{ord}$, $\varphi(x,\ov{a})$ defines a finite union of intervals and points with uniform bounds $K$ for the number of intervals and $P$ for the number of points. Thus, by taking $k_\varphi=2K+1$ and $C=P+1>0$ we can take the elements $c_i$ to be the end points of the intervals (including the minimal and maximal element), and the set $E$ of measures can be the finite set of all possible vectors $\ov{\mu}\in \{0,1\}^{k_\varphi}$.
\eej

\brmk \label{rmk-exact} Note that in the previous example the error term $C|[u,v]|^{1/2}$ is not necessary. In the context of asymptotic classes, this kind of phenomenon has been studied and there is a definition of \emph{asymptotic exact classes} of finite structures. A systematic study of these classes of finite structures will appear in \cite{AMSW}, which is work in preparation due to  S. Anscombe, D. Macpherson, C. Steinhorn and D. Wolf.
\ermk

The following is a technical lemma that will be useful to exhibit more examples of ordered asymptotic classes.

\bl \label{lemmaboolean} Let $\mathcal{C}$ be a class of finite linearly ordered structures, and suppose that conditions (1) and (2) of Definition \ref{defoclasses} hold for every formula in a subcollection $\Psi=\{\psi_s(x;\ov{z}_s):s<\omega\}\subseteq L$ that is closed under negations and intersections. Assume also that for every formula $\varphi(x,\ov{y})$, there are $\psi_{s_1}(x,\ov{z}_{s_1}),\ldots,\psi_{s_\ell}(x,\ov{z}_{s_k})\in\Psi$ such that for every $M\in\mathcal{C}$ \[M\models \forall \ov{y}\,\exists \ov{z}_{s_1},\ldots,\exists \ov{z}_{s_k}\,\forall x\left(\varphi(x,\ov{y})\leftrightarrow \bigvee_{i=1}^k \psi(x,\ov{z}_{s_i})\right).\] Then, $\mathcal{C}$ is an o-asymptotic class.
\el
\bdem Let $\varphi(x,\ov{y})$ be a formula in the language $L$, $M\in\mathcal{C}$ and $\ov{a}\in M^{|\ov{y}|}$. If $\varphi(x,\ov{y})\in \Psi$, we are done. Otherwise, since $\Psi$ is closed under negations and intersections, we can put $\ov{z}=\ov{z}_{s_1},\ldots,\ov{z}_{s_k}$ and replace the formulas $\psi_{s_i}$ by formulas $\psi'_{s_1}(x,\ov{z}):=\psi_{s_1}(x,\ov{z}_{s_1}), \psi'_{s_2}(x,\ov{z}):=\psi_{s_2}(x,\ov{z}_{s_1})\wedge \neg \psi'_{s_1}(x,\ov{z}),\ldots,\psi'_{s_k}(x,\ov{z}):=\psi_{s_k}(x,\ov{z}_{s_k})\wedge \neg \psi'_{s_1}(x,\ov{z})\wedge \cdots \wedge \neg \psi'_{s_{k-1}}(x,\ov{z}).$

By construction, for every $M\in\mathcal{C}$ and $\ov{a}\in M^{|\ov{y}|}$ there is a tuple $\ov{b}=\ov{b}_{s_1},\ldots,\ov{b}_{s_k}$ such that the sets $\psi'_{s_1}(M,\ov{b}),\ldots,\psi'_{s_\ell}(M,\ov{b})$ are disjoint and $\displaystyle{\varphi(M,\ov{a})=\bigcup_{i=1}^\ell \psi'_{s_i}(x,\ov{b})}$.

By hypothesis, each of the sets $\psi'_{s_i}(M,\ov{b})$ admits a decomposition given by elements \[c_0^{\psi'_{s_i}(x,\ov{b})}=\min M\leq c_1^{\psi'_{s_i}(x,\ov{b})}\leq \cdots \leq c_{k_{\psi'_{s_i}}}^{\psi'_{s_i}(x,\ov{b})}=\max M,\] with measures $\ov{\mu}^{\psi'_{s_i}(x,\ov{b})}$. By collecting all the elements $\langle c_r^{\psi'_{s_i}(x,\ov{b})}:i\leq \ell, r\leq k_{\psi'_{s_i}}\rangle$ and organizing them in increasing order, we find a decomposition of the structure $M$ into at most $k_\varphi:=k_{s_1}\cdots k_{s_\ell}$ intervals. Moreover, each interval $I$ in this decomposition is contained in some intersection of intervals $[c_{r_i-1}^{s_i},c_{r_i}^{s_i}]$ (for $i\leq \ell$), and the corresponding sum of measures $\mu_{r_1}^{\psi'_{s_1}(M,\ov{b})}+\cdots +\mu^{\psi'_{s_\ell}(M,\ov{b})}_{r_\ell}$ is at most $1$. So, we can take $E_\varphi\subseteq [0,1]^{k_\varphi}$ to be the finite set obtained by first adding zeros to all original tuples in the sets $E_{\psi'_{s_1}(x,\ov{z})},\ldots,E_{\psi'_{s_\ell}(x,\ov{z})}$ to make them of length $k_\varphi$ (in all possible places) and then take all the possible finite sums.\\

These choices are enough to ensure condition (1) of Definition \ref{defoclasses}. Finally, condition (2) follows by taking the conjunction of the corresponding formulas $\psi'_{s_i,\ov{\mu}}(\ov{y})$ witnessing the definability condition for $\psi'_{s_1}(x,\ov{y}),\ldots,\psi'_{s_\ell}(x,\ov{y})$.
\edem

\be \label{example-Ckcol} Consider the language $L=\{<,P_1,\ldots,P_k\}$ where $P_1,\ldots,P_k$ are unary predicates. Let $\mathcal{C}_{\text{$k$-col}}$ be the class of finite $L$-structures $M_n=([1,n],<,P_1,\ldots,P_k)$ where $M_n\models P_i(a)$ if and only if $a\equiv i \modulo{k}$, for each $n\geq k$.
\ee

\bp The class $\mathcal{C}_{\text{$k$-col}}$ is an o-asymptotic class.
\ep
\bdem By Lemma \ref{UQE-k-col}, every $L$-formula $\varphi(x,\ov{y})$ can be written as a finite disjunction of formulas of the form 
\[\varphi(x,\ov{y})=P_i(x)\wedge \bigwedge_{j\leq \ell}\tau_{1,j}(\ov{y})\leq x\leq \tau_{2,j}(\ov{y}).\] where $\tau_{1,j},\tau_{2,j}$ are terms in the language $L'=\{\min,\max,S,S^{-1}\}$. Thus, for every $M\in\mathcal{C}_{k-\text{col}}$ and $\ov{a}\in M^{|\ov{y}|}$, we can denote the interval $\bigcap_{j\leq \ell}[\tau_{1,j}(\ov{a}),\tau_{2,j}(\ov{a})]$ by $[d_1,d_2]$. Hence, the formula $\varphi(x,\ov{a})$ can be decomposed with measure $\ov{\mu}=(0,\frac{1}{k},0)$ into intervals by taking $c_0=\min M, c_1=d_1,c_2=d_2,c_3=\min M$, and the definability condition follows by choosing the end points of the interval (which depend only on $\ov{a}$). By Lemma \ref{lemmaboolean}, this is enough to show that $\mathcal{C}_{\text{$k$-col}}$ is an o-asymptotic class.
\edem

\section{Ultraproducts in o-asymptotic classes} \label{ultraoclasses}

In this section we present several results which will allow us to place infinite ultraproducts of structures in o-asymptotic classes in the classification theory map, i.e., we will study model theoretic properties such as o-minimality, NTP$_2$ and rosiness for infinite ultraproducts of o-asymptotic classes.

\bp \label{O-classes-o-minimal}Let $\mathcal{C}$ be a class of finite linearly ordered $L$-structures and suppose that every infinite ultraproduct is o-minimal. Then $\mathcal{C}$ is an o-asymptotic class.
\ep
\bdem We start with the following:

\noindent \emph{Claim:} For each formula $\varphi(x;\ov{y})$ there is $k<\omega$ such that for all $M_i\in\mathcal{C}$ and $\ov{a}\in M_i$,
\[M_i\models \neg \exists x_0,x_1,\ldots ,x_{2k+1}\left(\bigwedge_{1\leq i<j\leq 2k+1} (x_i<x_j) \wedge\bigwedge_{i=0}^k \varphi(x_{2i},\ov{a})\wedge \bigwedge_{i=0}^k \neg\varphi(x_{2i+1},\ov{a})\right).\]

\noindent \emph{Proof of the claim:} Assume otherwise. Then for every $k<\omega$ there are $M_{j_k}\in \mathcal{C}$, $\ov{a}^{j_k}\in M_{j_k}$ and elements $c_0^{j_k}<\cdots<c_{2k+1}^{j_k}$ such that
\[M_{j_k}\models \left(\bigwedge_{i=0}^k \varphi(c^{j_k}_{2i},\ov{a}^{j_k})\wedge \bigwedge_{i=0}^k \neg\varphi(c^{j_k}_{2i+1},\ov{a}^{j_k})\right).\] Let $\mathcal{U}$ be a non-principal ultrafilter on $\omega$ containing the set $J=\{j_k:k<\omega\}$. Then for $M=\prod_{\mathcal{U}} M_i$ and the tuples $\ov{a}=[\ov{a}^{j_k}]$, $c_i=[c_i^{j_k}]$ we have that $\{c_i:i<\omega\}$ is an increasing sequence in $M$ and $M\models \varphi(c_i,\ov{a})$ if and only if $i$ is even. Thus, the definable set $\varphi(M;\ov{a})$ cannot be expressed as a finite union of intervals and points, contradicting  o-minimality for $M$. \hspace{1cm}\checkmark\\

Let $\varphi(x,\ov{y})$ be an $L$-formula. From the claim above, it is clear that for every $M\in\mathcal{C}$ and $\ov{a}\in M$, the formula $\varphi(x;\ov{a})$ can be expressed as the union of at most $k=k_{\varphi}$ intervals and $k$ points, for a fixed $k<\omega$.

Define now the formula
\[\varphi(x,\ov{y},\ov{z}):=\varphi(x;\ov{y})\leftrightarrow \left(\bigvee_{i=0}^{k-1} z_i<x<z_{i+1} \vee \bigvee_{i=k+1}^{2k} x=z_i\right).\]
This formula holds in every structure in $\mathcal{C}$ (possibly repeating the intervals or the points) and therefore, for every infinite ultraproduct $M$ of structures in $\mathcal{C}$, we have
\[M\models \forall \ov{y}\, \exists z_0,\ldots,z_{2k}\forall x\left(\varphi(x;\ov{y},\ov{z})\right).\]

So, we can take $C=2k$, $E=\{\ov{\mu}\subseteq [0,1]^k: \mu_i\in \{0,1\}\}$ and for every $\ov{a}\in M$, we can take $\ov{c}=(c_0,\ldots,c_{k})$ where $\ov{c}'=(c_0,\ldots,c_k,c_{k+1},\ldots,c_{2k})$ is the corresponding tuple witnessing $M\models \forall x\left(\varphi(x;\ov{a},\ov{c})\right)$. This shows that the class $\mathcal{C}$ satisfies the condition (1) of the definition.
For the condition (2) (the definability clause), it is enough to take the formulas
$\varphi_{\ov{\mu}}(\ov{y};\ov{z}):=\forall x\left(\varphi(x,\ov{y},\ov{z})\right).$ \edem

\brmk As in Remark \ref{rmk-exact}, the previous result shows that if $\mathcal{C}$ is a class of linearly ordered structures whose infinite ultraproducts are o-minimal, then $\mathcal{C}$ is an o-asymptotic class and the error term $C|[u,v]|^{1/2}$ is not required because the approximation is exact.
\ermk
\brmk It is natural to ask whether Proposition \ref{O-classes-o-minimal} is true under conditions that are weaker than o-minimality, such as weak-o-minimality (every definable is a finite union of \emph{convex} sets and points) or quasi-o-minimality (cf. Definition \ref{def-quasi-o-min}). We have the following:
\benum
\item If $M$ is an ultraproduct of finite linearly ordered structures, then $M$ is discrete and every definable set in one variable has a minimal and a maximal element. Hence, $M$ is o-minimal if and only if it is weakly-o-minimal.
\item Proposition \ref{O-classes-o-minimal} is not longer true if we replace the condition of o-minimality by quasi-o-minimality. The class $\mathcal{C}_{PQ}$ defined in Section \ref{defCp} provides an example of a non o-asymptotic class whose ultraproducts are  quasi-o-minimal, as we showed in the Appendix \ref{Cpq-quasi-o-minimal}.
\eenum
\ermk

Since the leading idea in the definition of o-asymptotic classes is that they are melding properties from one-dimensional classes (whose ultraproducts are known to be simple and unstable in general) and o-minimal theories (which are known to be unstable theories with NIP), we are not expecting ultraproducts of o-asymptotic classes to be either simple or with NIP.  The two natural contexts which extend both simple and o-minimal theories are rosy theories and theories with NTP$_2$. We will show that both of these properties are satisfied by ultraproducts of o-asymptotic classes, starting with the following lemma which can be seen as an infinite version of condition (1) in Definition \ref{defoclasses}.

\bd Suppose that $M=\prod_{\mathcal{U}} M_n$ is an ultraproduct of finite structures, and let $D=\psi(M,\ov{a})$ be a non-empty definable set of $M$, with $\ov{a}=[\ov{a}_n]_\mathcal{U}$. We can define the \emph{counting measure on $M$ localized in $D$} to be the Keisler measure defined by 
\[\operatorname{meas}_D\left(\varphi(x,\ov{b})\right)=\lim_{n\to\mathcal{U}} \dfrac{|\varphi(M_n,\ov{b}_n)\cap \psi(M_n,\ov{a}_n)|}{|\psi(M_n,\ov{a}_n)|}.\] 
\ed
\bl \label{lemma-o-classes-measures} Let $\mathcal{C}$ be an o-asymptotic class, and let $M$ be an infinite ultraproduct of structures from $\mathcal{C}$. Then for every definable set $X\subseteq M^1$ there are a constant $C>0$, finitely many elements $c_0=\min M\leq c_1\leq \cdots \leq c_{k-1}\leq c_k=\max M$ and a tuple $\ov{\mu}=(\mu_1,\ldots,\mu_k)\in [0,1]^k$ such that for every $i\leq k$, either $\mu_i=0$ and $|X\cap [c_{i-1},c_i]|\leq C$ or for every infinite subinterval $[\alpha,\beta]\subseteq [c_{i-1},c_i]$ we have $\operatorname{meas}_{[\alpha,\beta]}(X)=\mu_i>0$.
\el
\bdem Suppose $M=\prod_\mathcal{U} M_n$ is an infinite ultraproduct of structures from $\mathcal{C}$ and $X$ is defined by $\varphi(x,\ov{b})$ for an $L$-formula $\varphi(x,\ov{y})$ and a tuple $\ov{b}=[\ov{b}_n]_{\mathcal{U}}\in M^{|\ov{y}|}$.   Since $\mathcal{C}$ is an o-asymptotic class, there is a constant $C=C_\varphi>0$ and an integer $k=k_\varphi$ and a finite set of tuples $E\subseteq [0,1]^k$ such that condition (1) in Definition \ref{defoclasses} holds, that is, for every $M_n$ there is a tuple $\ov{\mu}_n=(\mu_{1,n},\ldots,\mu_{k,n})\in E$ and elements $c_{0,n}=\min M_n\leq c_{1,n}\leq \cdots \leq c_{k-1,n}\leq c_{k,n}=\max M_n$ such that for every $i\leq k$, either $\mu_{i,n}=0$ or $\mu_{i,n}>0$ and for every $[u,v]\subseteq [c_{i-1},c_i]$, \[\left||\varphi(M_n;\ov{b}_n)\cap [u,v]|-\mu_i |[u,v]|\right|\leq C|[u,v]|^{1/2}.\]
Since $E$ is finite, there is a unique tuple $\ov{\mu}\in E$ such that $\ov{\mu}_n=\ov{\mu}=(\mu_1,\ldots,\mu_k)$ satisfies the conditions above for $\mathcal{U}$-almost all $n$. Let us consider the elements $c_0=[c_{0,n}]_\mathcal{U},\ldots,c_k=[c_{k,n}]_\mathcal{U}$ in $M$. 

If $\mu_i=0$, $M_n\models \exists^{\leq C}x\left(\varphi(x,\ov{b}_n)\wedge c_{i-1}\leq x\leq c_{i}\right)$, and by \Los' Theorem we will have $|X\cap [c_{i-1},c_i]|\leq C$. Otherwise, if $[\alpha,\beta]$ is an infinite interval of $[c_{i-1},c_i]$, then for $\mathcal{U}$-almost all $n$ we have
\begin{align*}
\mu_i-\dfrac{C}{|[\alpha_n,\beta_n]|^{1/2}}&\leq \dfrac{|\varphi(M_n,\ov{b})\cap [\alpha_n,\beta_n]|}{|[\alpha_n,\beta_n]|}\leq \mu_i+\dfrac{C}{|[\alpha_n,\beta_n]|^{1/2}}.
\end{align*}
Thus, since the interval $[\alpha,\beta]$ is infinite, we obtain 
$\operatorname{meas}_{[\alpha,\beta]}(X)=\operatorname{meas}_{[\alpha,\beta]}\left(\varphi(M;\ov{b})\right)=\mu_i>0$ by taking limits with respect to the ultrafilter $\mathcal{U}$.\edem


\bt \label{Uthrank1} Let $\mathcal{C}$ be an o-asymptotic class and let $M$ be an infinite ultraproduct of structures from $\mathcal{C}$. Then $\operatorname{Th}(M)$ is superrosy of $\uth$-rank $1$, i.e., the only formulas $\varphi(x,\ov{b})$ in $L(M)$ that \th-divide over the empty set are the algebraic formulas.
\et
\bdem Let $M$ be an infinite ultraproduct of structures in $\mathcal{C}$, and suppose $\varphi(x;\ov{b})$ is a non-algebraic formula that \th-divides over the empty set. Then there is a tuple of parameters $e\in M^{eq}$ such that $\operatorname{tp}(\ov{b}/e)$ is non-algebraic, and the set $\{\varphi(x,\ov{b}'):\ov{b}'\models \operatorname{tp}(\ov{b}/e)\}$ is $k$-inconsistent for some $k<\omega$. \\

Consider the finitely many elements $c^{\ov{b}}_0=\min M\leq c^{\ov{b}}_1\leq \cdots \leq c^{\ov{b}}_{\ell-1}\leq c^{\ov{b}}_\ell$ and the tuple $\ov{\mu}\in [0,1]^\ell$ ensured by Lemma \ref{lemma-o-classes-measures}. Using the definability condition, we may assume that $\operatorname{tp}(\ov{b}/e)\models \exists^{> C\cdot \ell}x(\varphi(x,\ov{y}))\wedge \exists z_0,\ldots,z_\ell(\varphi_{\ov{\mu}}(\ov{y};z_0,\ldots,z_\ell))$, where the first part of the disjunction ensures that $\varphi(x,\ov{b})$ is not algebraic while the second part states that $\varphi(x,\ov{b})$ admits a decomposition with proportion $\ov{\mu}$. Moreover, using the $\emptyset$-definable functions $m_{i,\ov{\mu}}^\varphi$ described in Remark  \ref{m-functions} we may suppose that $c_i^{\ov{b}}=m_{i,\ov{\mu}}^{\varphi}(\ov{b})$ for all $i=0,\ldots,\ell$.\\

Since $\phi(x,\ov{b})$ is not algebraic, $\ov{\mu}\neq (0,\ldots,0)$, and we can take $\mu_j$ to be the first non-zero coordinate of the tuple $\ov{\mu}$ such that $[c_{j-1}^{\ov{b}},c_j^{\ov{b}}]$ is infinite. To ease the notation, let us denote the functions $m_{j-1,\ov{\mu}}^{\varphi},m_{j,\ov{\mu}}^{\varphi}$ by $f_{j-1},f_j$ respectively. Thus, if $\ov{b}'\models \operatorname{tp}(\ov{b}/e)$ then $\varphi(x,\ov{b}')$ also has a decomposition with proportion $\ov{\mu}$, and the interval $[c_{j-1}^{\ov{b}'},c_j^{\ov{b}'}]=[f_{j-1}(\ov{b}'),f_j(\ov{b}')]$ is infinite. \\

\noindent \emph{\textbf{Claim:} There is an infinite sequence $\langle \ov{b}_s:s<\omega\rangle$ of realizations of $\operatorname{tp}(\ov{b}/e)$ and an infinite interval $[\alpha,\beta]$ such that $[\alpha,\beta]\subseteq [f_{j-1}(\ov{b}_s),f_j(\ov{b}_s)]$ for all $s<\omega$}\\

\noindent \emph{Proof of the Claim:}
First suppose there is $\alpha\in M$ such that $\alpha=f_{j-1}(\ov{b}_s)$ for countably many realizations $\ov{b}_s\models \operatorname{tp}(\ov{b}/e)$. Then, we can consider the type 
\[p(z):= \{S^n(\alpha)<w:n<\omega\}\cup\{ w< f_j(\ov{b}_s): s<\omega\}.\]
Note that this type is finitely satisfiable because the intervals $[f_{j-1}(\ov{b}_s),f_j(\ov{b})_s]$ are all infinite. Thus, by $\aleph_1$-saturation of $M$, there is a realization $\beta$ of $p(z)$, and we would have that $[\alpha,\beta]$ is an infinite interval contained in $[f_{j-1}(\ov{b}_s),f_j(\ov{b}_s)]$ for every $s<\omega$.\\

On the other hand, if for every $c\in M$ we have $c=f_{j-1}(\ov{b}')$ only for finitely many different realizations $\ov{b}'$ of $\operatorname{tp}(\ov{b}/e)$. Then, by compactness and $\aleph_1$-saturation, there is a uniform bound $t\in\mathbb{N}$ on the number of such realizations. Given a fixed formula $\theta(\ov{y})\in \operatorname{tp}(\ov{b}/e)$ we can define $F_\theta(z)=\exists \ov{y}\left(\theta(\ov{y})\wedge f_{j-1}(\ov{y})=z\right)$. The set $F_\theta(M) \supseteq \{f_{j-1}(\ov{b}'):\ov{b}'\models \operatorname{tp}(\ov{b}/e)\}$ is an infinite definable set, and has a decomposition into intervals given by elements $\min M=d_0^\theta<\cdots <d^{\theta}_{r}=\max M$ and a tuple of measures $\ov{\nu}^\theta=(\nu_1^\theta,\ldots,\nu_r^\theta)$.\\

Let $h\leq r$ be such that $f_{j-1}(\ov{b})\in [d^{\theta}_{h-1},d^{\theta}_{h}]$. Note that the elements $d^\theta_{h-1},d^\theta_h$ are the images of $\ov{b}$ under certain $e$-definable functions $g_{h-1},g_h$, because $f_{j-1}=m_{j-1,\ov{mu}}^\varphi$ is $\empty$-definable and the formula $F_\theta$ has parameters in $e$. So, since $d^{\theta}_{h-1}\leq f_{j-1}(\ov{y}\leq d^\theta_h$ is a formula in $\operatorname{tp}(\ov{b}/e)$, we will have $f_{j-1}(\ov{b}')\in [d^{\theta}_{h-1},d^{\theta}_{h}]$ for every $\ov{b}'\models \operatorname{tp}(\ov{b}/e)$.\\

We claim that $\nu_h^{\theta}>0$ and $[d_{h_1}^\theta,d_h^\theta]$ is infinite. Otherwise, $|F_\theta(M)\cap [d_{h-1}^\theta,d_h^\theta]|=C_{F_\theta}$ for some constant $C_{F_\theta}>0$ and we would have 
\begin{align*}
M\models &\exists z_1,\ldots,z_{C_{F_\theta}}\Bigg(\bigwedge_{i=1}^{C_{F_\theta}} (m^\theta_{h-1}(e)<z_i<m^\theta_{h}(e))\\
&\wedge \forall \ov{y}\Bigg(\theta(\ov{y})\wedge (m_{h-1}^\theta(e)<f_j(\ov{y})<m^\theta_{h}(e))\rightarrow \bigvee_{i=1}^{C_{F_\theta}} f_j(\ov{y})=z_i\Bigg)\Bigg)
\end{align*}
yielding that $\operatorname{tp}(\ov{b}/e)$ has at most $C_{F_\theta}\cdot t$ realizations,  which is absurd because $\operatorname{tp}(\ov{b}/e)$ is algebraic. Similarly, if $[d_{h-1}^\theta,d_h^\theta]$ were infinite, then $\operatorname{tp}(\ov{b}/e)$ would have at most $t\cdot |[d_{h-1}^\theta,d_h^\theta]|$ realizations, obtaining a contradiction.\\

We may assume without loss of generality that the interval $[f_{j-1}(\ov{b}),d_h^\theta]$ is infinite, because if $S^n(f_{j-1}(\ov{b})=d_h^\theta=g_h(e)$ then the same would be true for every $\ov{b}'\models \operatorname{tp}(\ov{b}/e)$, contradicting that every element $c\in M$ has at most $t$ preimages under the function $f_{j-1}$ that are realizations of $\operatorname{tp}(\ov{b}/e)$.\\

By the uniform distribution of $F_\theta(z)$ in the structures $M_n$ and \Los' Theorem we have that \[M\models \forall z\left(F_\theta(z) \wedge d^{\theta}_{h-1}<z<d^\theta_h \rightarrow \exists w (F_\theta(w)\wedge z<w<S^{D}(z))\right) \hspace{1cm} (**)\] where $D=\left(\dfrac{2C_{F_\theta}}{\nu_h}\right)^2$. Consider the type over $\ov{b},e$ given by \[p_1(\ov{y})=\{f_{j-1}(\ov{b})<f_j(\ov{y}) \}\cup \{\theta(\ov{y})\}\cup \{S^n(f_j(\ov{y}))<d^{\theta}_h, f_j(\ov{b}):n<\omega, \theta\in \operatorname{tp}(\ov{b}/e)\}.\]

This type is finitely satisfiable because given $\theta\in \operatorname{tp}(\ov{b}/e)$ and $n<\omega$, we can take the element $w$ provided by (**) when $z=f_{j-1}(\ov{b})$. Let $\ov{b}_1\models p_1$. Inductively, we may take $\ov{b}_{s+1}$ to be a realization of the type
\begin{align*}
p_{s+1}(\ov{y})=&\{f_{j-1}(\ov{b}_s)<f_{j-1}(\ov{y}) \}\cup \{\theta(\ov{y})\}\\
&\cup \{S^n(f_j(\ov{y}))<d^{\theta}_h, f_j(\ov{b}),f_j(\ov{b}_1),\ldots,f_j(\ov{b}_s):n<\omega, \theta\in \operatorname{tp}(\ov{b}/e)\}.
\end{align*}
It is easy to show now that the type
\[q(w,z)=\{f_{j-1}(\ov{b}_s)<w<z<f_j(\ov{b}_s):s<\omega\}\cup \{S^n(w)<z:n<\omega\}\] is finitely consistent, and by taking $(\alpha,\beta)\models q(w,z)$, we finish the proof of the Claim. \checkmark \\

Now, by Lemma \ref{lemma-o-classes-measures}, we have $\operatorname{meas}_{[\alpha,\beta]}(\varphi(x,\ov{b}_s))=\mu_j=\epsilon>0$ for every $s<\omega$, and by Fact \ref{kintersections}, there are $\ov{b}_{s_1},\ldots,\ov{b}_{s_k}$ with $s_1<\cdots<s_k$ such that
\[\operatorname{meas}_{[\alpha,\beta]}\left(\bigcap_{i=1}^k \varphi(x,\ov{b}_{s_i})\right)\geq \epsilon^{3^{k-1}}.\] In particular, $\{\varphi(x,\ov{b}_{s}):s<\omega\}$ is not $k$-inconsistent, and neither is the set $\{\varphi(x,\ov{b}'): \ov{b}'\models \operatorname{tp}(\ov{b}/e)\}$. This contradicts that $\varphi(x,\ov{b})$ strongly $k$-divides over $e$. 

\edem

\bt \label{O-NTP2} Every infinite ultraproduct of members of an o-asymptotic class is inp-minimal, and thus also NTP$_2$.
\et
\bdem Let $M$ be an infinite ultraproduct of structures in an o-asymptotic class, and suppose for a contradiction, using Fact \ref{mutuallyindisc}, that there are formulas $\varphi(x,\ov{y}),\psi(x,\ov{z})$ with $|x|=1$ and mutually indiscernible sequences $\langle \ov{a}_i:i<\omega\rangle,\langle \ov{b}_i:i<\omega\rangle$ witnessing an inp-pattern of depth $2$. That is:

\bitem
\item[(a)] For some $k$, both sets $\{\varphi(x;\ov{a}_i):i<\omega\}$ and $\{\psi(x;\ov{b}_i):i<\omega\}$ are $k$-inconsistent.
\item[(b)] For every $i,j<\omega$, $\varphi(M;\ov{a}_i)\cap \psi(M;\ov{b}_i)\neq \emptyset$.
\eitem

Using the definability conditions and indiscernibility, we may assume that there are integers $k_\varphi,k_\psi$ and tuples $\ov{\mu}\in [0,1]^{k_\varphi},\ov{\nu}\in [0,1]^{k_\psi}$ such that $\varphi(M;\ov{a}_i)$ admits a decomposition with proportion $\ov{\mu}$ for each $i<\omega$, and $\psi(M;\ov{b}_j)$ admits a decomposition with proportion $\ov{\nu}$ for each $j<\omega$. Furthermore, as in the previous proof, we may assume that for every $i<\omega$ the decomposition of $\varphi(M;\ov{a}_i)$ is given by the images of definable functions $f_r(\ov{a}_i):=m_{\varphi,r}^{\ov{\mu}}(\ov{a}_i)$ for $r\leq k_\varphi$. Likewise, for every $j<\omega$ the decomposition of $\psi(M;\ov{b}_j)$ is given by elements $g_s(\ov{b}_j)$, with $s\leq k_\psi$.\\

By property (b), the definable set $\varphi(M,\ov{a}_i)\cap \psi(M,\ov{b}_j)$ is non-empty, and so it has a minimal element $\alpha_{i,j}$. Let $r\leq k_\varphi$ and $s\leq k_\psi$ be such that $\alpha_{i,j}\in [f_{r-1}(\ov{a}_i),f_r(\ov{a}_i)]\cap [g_{s-1}(\ov{b}_j),g_s(\ov{b}_j)]$. Note that both $r$ and $s$ are fixed by mutual indiscernibility and the definability of the functions $f_r,g_s$.\\

\emph{Claim 1:} We have $\mu_r,\nu_s>0$, and for every $i,j<\omega$ the intersection $[f_{r-1}(\ov{a}_i),f_r(\ov{a}_i)]\cap [g_{s-1}(\ov{b}_j),g_s(\ov{b}_j)]$ is infinite.\\

\emph{Proof of Claim 1:} If $\mu_r=0$, then $|\varphi(M,\ov{a}_i)\cap [f_{r-1}(\ov{a}_i),f_r(\ov{a}_i)]|\leq C_\varphi$, and for every $j<\omega$ we have $\alpha_{i,j}\in\psi(M,\ov{b}_j)\cap \varphi(M,\ov{a}_i)\cap [f_{r-1}(\ov{a}_i),f_r(\ov{a}_i)]$. Thus, by the pigeonhole principle, there is a single element $\alpha_i$ such that $\alpha_i=\alpha_{i,j}$ for infinitely many $j<\omega$, contradicting that $\{\psi(x,\ov{b}_j):j<\omega\}$ is $k$-inconsistent. We can use the same argument if the intersection $[f_{r-1}(\ov{a}_i),f_r(\ov{a}_i)]\cap [g_{s-1}(\ov{b}_j),g_s(\ov{b}_j)]$ is finite. \checkmark\\

\emph{Claim 2:} At least one of the following properties holds:
\benum
\item[(i)] For every $i_1,i_2<\omega$, the intersection $[f_{r-1}(\ov{a}_{i_1}),f_r(\ov{a}_{i_1})] \cap [f_{r-1}(\ov{a}_{i_2}),f_r(\ov{a}_{i_2})]$ is infinite.
\item[(ii)] For every $j_1,j_2<\omega$, the intersection $[g_{s-1}(\ov{b}_{j_1}),g_s(\ov{b}_{j_1})] \cap [g_{s-1}(\ov{b}_{j_2}),g_s(\ov{b}_{j_2})]$ is infinite.
\eenum
\emph{Proof of Claim 2:} First, notice that if $|[f_{r-1}(\ov{a}_{1}),f_r(\ov{a}_{1})] \cap [f_{r-1}(\ov{a}_{2}),f_r(\ov{a}_{2})]|=\ell$ for some integer $\ell\geq 1$, then the intersection $[f_{r-1}(\ov{a}_{1}),f_r(\ov{a}_{1})] \cap [f_{r-1}(\ov{a}_{2}),f_r(\ov{a}_{2})]$ is precisely the last $\ell$ elements of the interval $[f_{r-1}(\ov{a}_{1}),f_r(\ov{a}_{1})]$. By indiscernibility, the same is true for the intersection $[f_{r-1}(\ov{a}_{1}),f_r(\ov{a}_{1})] \cap [f_{r-1}(\ov{a}_{3}),f_r(\ov{a}_{3})]$, but then the intersection $[f_{r-1}(\ov{a}_2,f_r(\ov{a}_2)]\cap [f_{r-1}(\ov{a}_{3}),f_r(\ov{a}_{3})]$ contains the first $\ell$ elements of $[f_{r-1}(\ov{a}_2,f_r(\ov{a}_2)]$, contradicting that $\operatorname{tp}(a_1,a_3)=\operatorname{tp}(a_2,a_3)$. A similar argument shows that the intersection $[g_{s-1}(\ov{b}_1),g_s(\ov{b}_1)]\cap [g_{s-1}(\ov{b}_2),g_s(\ov{b}_2)]$ is not finite either.\\

Now suppose that both $[f_{r-1}(\ov{a}_{1}),f_r(\ov{a}_{1})] \cap [f_{r-1}(\ov{a}_{2}),f_r(\ov{a}_{2})]$ and $[g_{s-1}(\ov{b}_1),g_{s}(\ov{b}_1)]\cap [g_{s-1}(\ov{b}_2),g_s(\ov{b}_2)]$ are empty. As a first case, assume that $f_r(\ov{a}_1)<f_{r-1}(\ov{a}_2)$ and $g_{s}(\ov{b}_1)<g_{s-1}(\ov{b}_2)$. Then, since $\alpha_{12}$ belongs to $[f_{r-1}(\ov{a}_{1}),f_r(\ov{a}_{1})] \cap [g_{s-1}(\ov{b}_2),g_s(\ov{b}_2)]$, we have $g_{s}(\ov{b}_1)<g_{s-1}(\ov{b}_2)\leq f_r(\ov{a}_1)<f_r(\ov{a}_2)$, and so the intersection $[f_{r-1}(\ov{a}_2),f_r(\ov{a}_2)]\cap [g_{s-1}(\ov{b}_1),g_s(\ov{b}_1)]$ is empty, contradicting the fact that $\alpha_{21}$ belongs to this intersection. The other cases can be analyzed in a similar fashion, using a possibly different choice of the elements $\alpha_{11},\alpha_{12},\alpha_{21},\alpha_{22}$.\\

Therefore, the only possibility is that either $[f_{r-1}(\ov{a}_1),f_r(\ov{a}_1)]\cap [f_{r-1}(\ov{a}_2),f_r(\ov{a}_2)]$ is infinite or  $[g_{s-1}(\ov{b}_{1}),g_s(\ov{b}_{1})] \cap [g_{s-1}(\ov{b}_{2}),g_s(\ov{b}_{2})]$ is infinite. By indiscernibility, this implies that at least one of the properties (i) or (ii) holds.
\checkmark\\

Suppose without loss of generality that property (i) of Claim 2 holds. Note that for every $R<\omega$, the intersection \[[g_{s-1}(\ov{b}_1),g_s(\ov{b}_1)]\cap \bigcap_{i=1}^R [f_{r-1}(\ov{a}_i),f_r(\ov{a}_i)]\]
is equal to the intersection of two of these intervals, which by Claims 1 and 2 is infinite. Therefore, by compactness and $\aleph_1$-saturation, we can find an infinite interval $[\alpha,\beta]\subseteq [f_{r-1}(\ov{a}_i),f_r(\ov{a}_i)]\cap [g_{s-1}(\ov{b}_1),g_s(\ov{b}_1)]$ for all $i<\omega$. Hence, by Lemma \ref{lemma-o-classes-measures} we have $\operatorname{meas}_{[\alpha,\beta]}(\varphi(x,\ov{a}_i))=\mu_r>0$ for all $i<\omega$, and by Fact \ref{kintersections} we would have \[\operatorname{meas}_{[\alpha,\beta]}\left(\bigcap_{s=1}^k\varphi(x,\ov{a}_{i_s})\right)\geq \mu_r^{3^{k-1}}>0\]for some $i_1<\ldots<i_k<\omega$, contradicting $k$-inconsistency of $\{\varphi(x,\ov{a}_i):i<\omega\}$ in (a).\\

If property (ii) of Claims 2 holds, then we can use a symmetric argument and contradict $k$-inconsistency of $\{\psi(x,\ov{b}_j):j<\omega\}$. 
\edem
\section{Cyclic groups with an ordering} \label{examples-O-classes}

It was already mentioned in Section \ref{odefcelldecomp} that the class of finite linear orders is an o-asymptotic class. We will present in this section an example of an o-asymptotic class with some algebraic features.

\bd Given a natural number $N$, we consider the finite linearly ordered structure $\mathcal{Z}_N=(\mathbb{Z}/(2N+1)\mathbb{Z},\ci{+},\ci{<})$ where $\ci{+}$ denotes the usual addition in the cyclic group and the linear order $\ci{<}$ is imposed as:
\[-\ov{N}\ \ci{<}\ -\ov{(N-1)}\ \ci{<}\cdots \ci{<}\ \ov{0}\ \ci{<}\cdots \ci{<}\ \ov{N-1} \ \ci{<}\ \ov{N},\]
where $\ov{k}$ denotes the equivalence class of $k$ modulo $2N+1$.
\ed
During this section we will show that the class $\mathcal{C}_{ocyc}=\{\mathcal{Z}_N:N<\omega\}$ is an o-asymptotic class. For this, we first need a result describing uniformly the definable sets for structures in this class. The general idea will be to use the natural interpretation of the structures $\mathcal{Z}_N$ into $(\mathbb{Z},+,-,0,1,<)$ and use Presburger's Theorem to pull back a uniform description for the definable sets in the structures $\mathcal{Z}_N$. We now proceed to describe this idea in more detail.\\ 

\subsection{Quantifier elimination results in Presburger arithmetic} To start, let us recall the following well-known result on quantifier elimination for the structure $(\mathbb{Z},+,-,<,0,1)$.
\bh [Presburger's Theorem, as presented in {\cite{BPW}}] Every formula $\varphi(x;\ov{y})$ in the language $L=\{0,1,+,-,<\}$ is equivalent in $(\mathbb{Z},+,-,0,1<)$ to a boolean combination of formulas of the form $n\cdot x=t(\ov{y})$, $n\cdot x<t(\ov{y})$, $n\cdot x>t(\ov{y})$ or $D_m(n\cdot x+t(\ov{y}))$, where $t(\ov{y})$ are terms depending only on $\ov{y}$ and $D_m(z)$ is the formula $\exists t(\underbrace{t+\cdots +t}_{\text{$m$ times}}=z)$.
\eh
Regarding the formulas $D_m(n\cdot x+y)$, we have the following simplifications, which are fairly standard when studying definable sets in Presburger Arithmetic.

\bh \label{simplif1} For every $m\geq 2$, $n,b\in\mathbb{Z}$ there are $m'\geq 2$ and $0\leq b'\leq m'-1$ such that the formulas $D_m(n\cdot x+b)$ and $D_{m'}(x-b')$ define the same set in $(\mathbb{Z},+,-,<,0,1)$.
\eh
\bh \label{simplif2} For every $m_1,m_2\geq 2$ and integers $0\leq b_1<m_1,0\leq b_2<m_2$, either the conjunction $D_{m_1}(x-b_1)\wedge D_{m_2}(x-b_2)$ defines the empty set or there are $m\geq 2$ and $0\leq b<m$ such that $D_{m_1}(x-b_1)\wedge D_{m_2}(x-b_2)$ defines the same set as the formula $D_m(x-b)$ in $(\mathbb{Z},+,-,<,0,1)$.
\eh

\subsection{Interpretation of the structures in $\mathcal{C}_{ocyc}$ in Presburger arithmetic}\ \\

In this subsection, we will analyze the definable sets along the class $\mathcal{C}_{ocyc}$ and use interpretations in $(\mathbb{Z},+,-,0,1,<)$ to obtain a uniform quantifier elimination result for $\mathcal{C}_{ocyc}$ (see Proposition \ref{qeocyc}). 

Let us consider the language $L'=\{0,\ci{+},\ci{<},\min,\max\}$ and for every $N\geq 1$, the function $f:\mathcal{Z}_N\to\mathbb{Z}$ given by $f(\overline{x})=x$ for every $x\in [-N,N]$. We can interpret each symbol of $L'$ in the structure $(\mathbb{Z},+,-,0,1,<)$ via the following:
\begin{align*}
&\text{\emph{Constants:} } \min:=-N, \max:=N.\\
&\text{\emph{Formulas:} } \ Z(x;N):=-N\leq x\leq N \equiv \min\leq x\leq \max .\\
&\text{\emph{Order:} } \hspace{0.75cm}O(x,y;N):=Z(x;N)\wedge Z(y;N)\wedge x<y.\\
&\text{\emph{Addition}:\ \ }\hspace{0.2cm}S(x,y,z;N):=``x\ \ci{+}\ y=z \text{\ \ in $\mathcal{Z}_N$}"\\
&:=Z(x;N)\wedge Z(y;N) \wedge Z(z;N) \wedge [x+y=z \vee x+y=z+(2N+1) \vee x+y=z-(2N+1)]\\
&\equiv Z(x;N)\wedge Z(y;N) \wedge Z(z;N)\wedge\big[(z<0<x,y \wedge (x+\min)+(y+\min)=z+1)\\
&\vee (\neg (z<0<x,y \wedge x,y<0<z)\wedge (x+y=z)\vee (x,y<0<z \wedge (x+\max)+(y+\max)+1=z)\big].\\
\end{align*}

This allows us to find a translation in $(\mathbb{Z},+,-,0,1,<)$ of the $L'$-formulas coming from the structures $\mathcal{Z}_N$. Namely, we can identify $\mathcal{Z}_N$ with its image $f(\mathcal{Z}_N)=[-N,N]\subseteq \mathbb{Z}$, and for every $L'$-formula $\varphi(x_1,\ldots,x_n)$ there is a formula $\widehat{\varphi}(x_1,\ldots,x_n;w)$ in the language $L=\{+,<,-,0\}$ such that $\mathcal{Z}_N\models \varphi(a_1,\ldots,a_n)\text{ if and only if }(\mathbb{Z},+,-,0,1,<)\models \widehat{\varphi}(f(a_1),\ldots,f(a_n),N)$. Furthermore, for every choice of parameters $\ov{b},N$, we have $\varphi(\mathcal{Z}_N;\ov{b})=[-N,N]\cap\widehat{\varphi}(\mathbb{Z};\ov{b},N)$.\\

Given an $L'$-formula of the form $\varphi(x;\ov{y})$, we know by Presburger's Theorem and Lemmas \ref{simplif1}, \ref{simplif2} that the formula $\widehat{\varphi}(x,\ov{y},w)$ is a boolean combination of formulas of the form $n\cdot x=t(\ov{y},w)$, $n\cdot x<t(\ov{y},w)$, $n\cdot x>t(\ov{y},w)$ or $D_m(x+t(\ov{y},w))$, where $t(\ov{y},w)$ are terms depending only on $\ov{y},w$ and $D_m(z)$ is the formula $\exists t(\underbrace{t+\cdots +t}_{\text{\tiny{$m$ times}}}=z)$. 

Thus, in order to show that $\mathcal{C}_{ocyc}$ is an o-asymptotic class, it is necessary to show that for every choice of parameters $\ov{b}\in \mathcal{Z}_N=[-N,N]\subseteq \mathbb{Z}$, the finite intersections of the formulas $n\cdot x=t(\ov{b},N)$, $n\cdot x<t(\ov{b},N)$, $n\cdot x>t(\ov{b},N)$, $D_m(x+t(\ov{b},N))$ satisfy three conditions: they can be defined using the language $L'$, they define sets that can be decomposed into a uniformly bounded number of intervals (i.e., the number of intervals does not depend on $N$), and finally, the extreme points of such intervals can be defined from $\ov{b}$.\\

To show these conditions, we will work in the structure $(\mathbb{Z},+,-,<,0,1)$, using the identification $\mathcal{Z}_N=[-N,N]\subseteq \mathbb{Z}$. To clarify the notation,  we write $\ci{+},\ci{<}$ to distinguish the operations and relations in $\mathcal{Z}_N$ from the corresponding operations and relations in $\mathbb{Z}$. Also, to emphasize the uniformity of these interpretations, let us use the variable $w$ instead of the integer $N$.

\brmk \label{oplusvssum} For $a,b\in \mathcal{Z}_w$ both $a\ \ci{+} b$ and $a+b$ are well defined integers, although $a+b$ might not be an element in $\mathcal{Z}_w$. In fact, by modular addition, we have the following identities:
\[\begin{cases}a+b=(a\oplus b)-(2w+1), &\text{if $a+b<a\oplus b$, or equivalently, if $a\oplus b>a$ and $b<0$}\\
a+b=(a\oplus b)+(2w+1), &\text{if $a+b>a\oplus b$, or equivalently, if $a\oplus b<a$ and $b>0$}\\a+b=a\oplus b, &\text{otherwise}.\end{cases}\]
\ermk

\bd Define recursively the \emph{carry} of a tuple $\ov{y}=(y_1,\ldots,y_n)\in(\mathcal{Z}_w)^n$ as follows:
\bitem
\item $\operatorname{carry}(y_1)=0.$
\item $\operatorname{carry}(y_1,\ldots,y_{k+1})=\begin{cases}\operatorname{carry}(y_1,\ldots,y_k)+1,&\text{if $\displaystyle{\bigoplus_{i=1}^{k+1}y_i<\bigoplus_{i=1}^k y_i}$ and $y_{k+1}>0$}\\ \operatorname{carry}(y_1,\ldots,y_k)-1,&\text{if $\displaystyle{\bigoplus_{i=1}^{k+1}y_i>\bigoplus_{i=1}^k y_i}$ and $y_{k+1}<0$}\\ \operatorname{carry}(y_1,\ldots,y_k),&\text{otherwise}.\end{cases}$
\eitem
\ed

This ``carry'' will be an important tool to translate definable sets from $(\mathbb{Z},+,-,0,1,<)$ to the structure $(\mathcal{Z}_w,\min,\max,\ci{+},\ci{<},0)$. Intuitively, it represents the number of times we have to ``wrap around'' in the cyclic group $\mathcal{Z}_w$ while performing the addition $y_1\ \ci{+}\ldots\ci{+}\ y_n$. From this point of view, the following lemma should be clear:
\bl Given $y_1,\ldots,y_n\in\mathcal{Z}_w$, we have $\displaystyle{\sum_{i=1}^n y_i=\left(\operatorname{carry}(y_1,\ldots,y_n)\right)\cdot (2w+1) + \bigoplus_{i=1}^n y_i.}$
\el

\bdem The statement follows easily by induction on $n$.
\edem

Note that for every $n,k$, the condition ``$\operatorname{carry}(y_1,\ldots,y_n)=k$'' is uniformly $\emptyset$-definable in the class $\mathcal{C}_{ocyc}=\{\mathcal{Z}_w:1\leq w<\omega\}$. For instance, when $k\geq 0$ and $y_1,\ldots,y_n\geq 0$, it is defined by
\[\bigvee_{I\subseteq \{1,\ldots,n\},|I|=k}\left(\bigwedge_{i\in I}\bigoplus_{j=1}^{i+1}y_j<\bigoplus_{j=1}^{i}y_j \wedge \bigwedge_{i\not\in I}\bigoplus_{j=1}^{i}y_j\leq \bigoplus_{j=1}^{i+1}y_j\right).\]

On the other hand, the formula $\operatorname{carry}(\underbrace{x,\ldots,x}_{\text{$n$ times}})=k$ defines an interval in $\mathcal{Z}_w$. In fact, we have 
\begin{align*}
\operatorname{carry}(\underbrace{x,\ldots,x}_{n})=k&\Leftrightarrow -w+k\cdot (2w+1)\leq n\cdot x\leq w+k(2w+1)\\
&\Leftrightarrow \dfrac{-w+k\cdot (2w+1)}{n}\leq x\leq \dfrac{w+k(2k+1)}{n}\\
&\Leftrightarrow x\in \left[\dfrac{-w+k\cdot (2w+1)}{n},\dfrac{w+k\cdot (2w+1)}{n}\right]\cap [-w,w].
\end{align*}
The end points of this interval are also definable. For example, the maximal element is given by the formula 
\[\theta_{k,+}(x):=\operatorname{carry}(\underbrace{x,\ldots,x}_n)=k \wedge \left(x=\max \ \vee \operatorname{carry}(\underbrace{S(x),\ldots,S(x)}_n)\neq k\right).\]

Now, we start analyzing each of the possible atomic formulas appearing in the boolean decomposition of $\widehat{\varphi}(x;\ov{y},w)$, starting with $nx<t(\ov{y},w)$. First, notice that we can assume (possibly changing the signs of some variables in the tuple $(\ov{y},w))$ that $t(\ov{y},w)=\sum_{i=1}^{|\ov{y}|} \alpha_i\cdot y_i + \beta\cdot w$, with $\alpha_i,\beta\in \mathbb{Z}^{\geq 0}$. To ease the notation, let us assume that $t(\ov{y},w)=\alpha\cdot y+\beta\cdot w$ for a single variable $y$, and denote by $t_\oplus(y,w)$ the term obtaining by replacing every appearence of the function symbol $+$ by $\ci{+}$. Then,

\begin{align*}
&\mathbb{Z}\models n\cdot x< \alpha\cdot y +\beta\cdot w \Leftrightarrow \mathbb{Z}\models \underbrace{x+\cdots+x}_{n}<(\underbrace{y+\cdots+y}_{\alpha})+(\underbrace{w+\cdots+w}_{\beta})\\
\Leftrightarrow \ &\underbrace{\bigoplus_{i=1}^n x}_{\text{in $\mathcal{Z}_w$}} + \operatorname{carry}(\underbrace{x,\ldots,x}_n)\cdot (2w+1)<\operatorname{carry}(\underbrace{y,\ldots,y}_\alpha;\underbrace{w,\ldots,w}_\beta)\cdot (2w+1)+\underbrace{\bigoplus_{\alpha}y\oplus\bigoplus_{\beta}w}_{\text{in $\mathcal{Z}_w$}}.\\
\end{align*}
So, either $\operatorname{carry}(\underbrace{x,\ldots,x}_n)<\operatorname{carry}(\underbrace{y,\ldots,y}_\alpha,\underbrace{w,\ldots,w}_\beta)=\operatorname{carry}(\underbrace{y,\ldots,y}_\alpha,\underbrace{\max,\ldots,\max}_\beta)$, or they are equal and $\displaystyle{\bigoplus_{i=1}^n x<t_{\oplus}(y,w)=t_{\oplus}(y,\max)}$. Hence,

\begin{align*}
&\mathbb{Z}\models Z(x,w)\wedge Z(y,w)\wedge n\cdot x<\alpha\cdot y+\beta\cdot w\\
\Leftrightarrow \hspace{0.5cm} &\mathcal{Z}_w\models \bigvee_{k_1<k_2\leq \max\{n,\alpha+\beta\}} \left[\operatorname{carry}(\underbrace{x,\ldots,x}_n)=k_1 \wedge \operatorname{carry}(\underbrace{y,\ldots,y}_\alpha,\underbrace{\max,\ldots,\max}_\beta)=k_2\right] \\
&\hspace{-1cm}\vee \bigvee_{-n\leq k\leq n} \left[\operatorname{carry}(\underbrace{x,\ldots,x}_n)=k=\operatorname{carry}(\underbrace{y,\ldots,y}_\alpha,\underbrace{\max,\ldots,\max}_\beta) \wedge \left(\underbrace{x\ \ci{+}\ldots\ci{+}\ x}_n<t_\oplus(y,\max)\right)\right].
\end{align*}

Note that the previous formula does not mention the variable $w$. So, for every $\ov{b}\in\mathcal{Z}_w\subseteq \mathbb{Z}$, the formula $n\cdot x<t(\ov{b},w)$ defines in $\mathcal{Z}_w$ a union of at most $2n+\max\{n,\alpha+\beta\}$ disjoint intervals, whose end points are uniformly definable from $\ov{b}$ in the language $L'$. We can use a similar argument for the formulas $n\cdot x>t(\ov{y},w)$ and $n\cdot x=t(\ov{y},w)$. \\

Now we analyze the formula $D_m(x+t(\ov{y},w))$. Again, let us assume $t(\ov{y},w)=\alpha\cdot y+\beta\cdot w$. Note first that for the formula $D_m(x)$, we have that 
$\mathcal{Z}_N\cap D_m(\mathbb{Z})=P_m(\mathcal{Z}_N)$, where $P_m$ is the $L'$-formula 
\[P_m(x):=\exists t\left([0\,\text{\ci{<}}\,t\,\text{\ci{<}}\,t\,\text{\ci{+}}\,t\,\text{\ci{<}}\cdots\text{\ci{<}}\,\underbrace{t\,\text{\ci{+}}\cdots\text{\ci{+}}\,t}_m=x] \vee [\underbrace{t\,\text{\ci{+}}\cdots\text{\ci{+}}\,t}_m=x\,\text{\ci{<}}\cdots\text{\ci{<}}\,t\text{\ci{+}}\,t\,\text{\ci{<}}\,t\,\text{\ci{<}}\,0]\right).\]
Let $\alpha',\beta'$ be elements such that $\alpha\equiv \alpha',\beta\equiv\beta' \modulo{m}$. Also, note that there are unique elements $0\leq y',w'<m$ such that $\mathbb{Z}\models D_m(y-y')\wedge D_m(w-w')$, and we have
\begin{align*}
&\mathbb{Z}\models D_m(x+t(\ov{y},w))\Leftrightarrow \mathbb{Z}\models \bigvee_{y'=0}^{m-1}\bigvee_{w'=0}^{m-1} (D_m(y-y')\wedge D_m(w-w')\wedge D_m(x+\alpha'\cdot y'+\beta'\cdot w'))\\
&\Leftrightarrow \mathcal{Z}_w \models \bigvee_{y'=0}^{m-1}\bigvee_{w'=0}^{m-1} (P_m(y-y')\wedge P_m(w-w')\wedge P_m(x+\alpha'\cdot y'+\beta'\cdot w')).
\end{align*} 
We summarize these results in the following proposition.
\bp \label{qeocyc} Consider the language $L'=\{\text{\ci{+}},\text{\ci{-}},\text{\ci{<}},0,1,\min,\max\}$. Every $L'$-formula $\varphi(x,\ov{y})$ is equivalent (uniformly in $\mathcal{C}_{ocyc}$) to a boolean combination of formulas of the form $x=t(\ov{y})$, $x<t(\ov{y})$ and $P_m(x+t(\ov{y}))$ where $t(\ov{y})$ is a definable function in the language $L'$.
\ep
\bdem Given the $L'$-formula $\varphi(x,\ov{y})$, there is a formula $\widehat{\varphi}(x;\ov{y},w)$ in the language $L=\{+,-,<,0,1\}$ such that $\mathbb{Z}\models \widehat{\varphi}(x;\ov{y},w)$ iff $\mathcal{Z}_w\models \varphi(x;\ov{y})$. By Presburger's Theorem and Lemmas \ref{simplif1} and \ref{simplif2}, $\widehat{\varphi}(x;\ov{y},w)$ is equivalent in $\mathbb{Z}$ to a boolean combination of formulas of the form $n\cdot x<t(\ov{y},w), n\cdot x=t(\ov{y},w)$ and $D_m(x+t(\ov{y},w))$.\\

By the previous discussion in this section, for each of these formulas there is a boolean combination of $L'$-formulas of the form $x<t(\ov{y}),x=t(\ov{y})$ and $P_m(x+t(\ov{y}))$, where $t(\ov{y})$ is an $L'$-definable function.\edem

\bt The class $\C_{ocyc}$ is an o-asymptotic class.
\et
\bdem Note first that for every $N\geq 1$, $Z_N\models \neg P_m(z)\Leftrightarrow \bigvee_{i=1}^{m-1}P_m(z+i)$. Hence, by the previous lemma and Lemma \ref{lemmaboolean}, it is enough to check the conditions in Definition \ref{defoclasses} for finite conjunctions of the form \[\varphi(x,\ov{y}):=\bigwedge_i t_1^i(\ov{y})\leq x\leq t_2^i(\ov{y})\wedge \bigwedge_{j}P_{m_j}(x+t^j(\ov{y})),\] where $t_1^i(\ov{y}),t_2^i(\ov{y}),t^j(\ov{y})$ are $L'$-definable functions. Moreover, by Lemma \ref{simplif1}, the second conjunction is equivalent to a single formula $P_m(x+t(\ov{y}))$, while the first conjunction defines a single interval $[d_1(\ov{y}),d_2(\ov{y})]$. \\

Thus, we can take $k_\varphi=3$, $c_0=\min\leq c_1=d_1(\ov{y})\leq c_2=d_2(\ov{y})\leq c_3=\max$, and the tuple of measures $\ov{\mu}=(0,\frac{1}{m},0)$ (if the last conjunction does not appear, then $\ov{\mu}=(0,1,0)$. Finally, the definability condition is given by the formula \[\varphi_{\ov{\mu}}(\ov{y};z_0,z_1,z_2,z_3):=(z_0=\min)\wedge (z_3=\max) \wedge (z_1=d_1(\ov{y})) \wedge (z_2=d_2(\ov{y})).\qedhere\] 
\edem

\section{Non-examples of o-asymptotic classes.} In this section we will describe some examples of classes of linearly ordered finite structures which are not o-asymptotic classes. Often, in order to show that a class $\mathcal{C}$ is not o-asymptotic, we will give a description of their infinite ultraproducts and use one of the following two results.

\bl \label{lemmaalternateintervals} Suppose that $\mathcal{C}$ is a class of finite linearly ordered structures. Suppose that in some ultraproduct $M=\prod_\mathcal{U} M_n$ of structures in $\mathcal{C}$ there is a definable set $X\subseteq M^1$ and infinite convex sets $\langle I_s,J_s:s<\omega\rangle$ such that for every $s<\omega$ we have $I_s<J_s<I_{s+1}$, $I_s\subseteq X$ and $J_s\subseteq M\setminus X$.  Then, $\mathcal{C}$ is not an o-asymptotic class.
\el
\bdem Suppose that $X=\varphi(M;\ov{b})$ for some formula $\varphi(x,\ov{y})$. By Lemma \ref{lemma-o-classes-measures} there are finitely many elements $c_0=\min M\leq c_1\leq \cdots \leq c_{k-1}\leq c_k=\max M$ and a tuple $\ov{\mu}\in [0,1]^k$ such that whenever $[c_{i-1},c_i]\cap X$ is infinite, $\operatorname{meas}_{[\alpha,\beta]}(X)=\mu_i$ for every infinite subinterval $[\alpha,\beta]\subseteq [c_{i-1},c_i]$. By the pigeonhole principle, one of the intervals $[c_{i-1},c_i]$ contains both $I_s,J_s$ for some $s<\omega$. We then have that $[c_{i-1},c_i]\cap X\supseteq I_s$ is infinite, and by compactness there is an infinite interval $[\alpha,\beta]\subseteq J_s$. Thus, we would have $\operatorname{meas}_{[\alpha,\beta]}(X)=\operatorname{meas}_{[\alpha,\beta]}(X\cap J_s)=0$, a contradiction.\edem

\bl \label{lemmameasurezero} Let $\mathcal{C}$ be a class of finite linearly ordered $L$-structures and $M$ an infinite ultraproduct of structures in $\mathcal{C}$. Suppose that there is an infinite definable set $X\subseteq M^1$ such that for every infinite interval $[\alpha,\beta]$, $\operatorname{meas}_{[\alpha,\beta]}(X)=0$. Then, $\mathcal{C}$ is not an o-asymptotic class.
\el
\bdem If $\mathcal{C}$ were an o-asymptotic class, then by Lemma \ref{lemma-o-classes-measures} there would be finitely many elements $c_0=\min M\leq c_1\leq \cdots \leq c_k=\max M$ and a tuple $\ov{\mu}=(\mu_1,\ldots,\mu_k)\in [0,1]^k$ such that for every $i\leq k$, either $\mu_i=0$ and $|X\cap [c_{i-1},c_i]|\leq C$, or $\mu_i>0$ and for every infinite $[u,v]\subseteq [c_{i-1},c_i]$, $\operatorname{meas}_{[u,v]}(X)=\mu_i$. Since $X$ is infinite, we must have $\mu_i\neq 0$ and $[c_{i-1},c_i]$ infinite for some $i\leq k$. Thus,  by taking $[\alpha,\beta]=[c_{i-1},c_i]$ we would have $\operatorname{meas}_{[\alpha,\beta]}(X)=\mu_i\neq 0$, contradicting the assumptions.
\edem

We start now describing non-examples, some of which might be expected to be o-asymptotic classes but are not. The following can be understood as a prototypical minimal counterexample.

\bej \label{defCp} Let $L$ be the language $L=\{<,P,Q\}$ where $P,Q$ are unary predicates. Consider the finite $L$-structures $M_n$ given by $M_n=([1,n\cdot 2n^3],<,P,Q)$ where $<$ is the usual order on $\Nat$ and the predicates $P$ and $Q$ are interpreted as follows:
\begin{align*}
P(M_n)&:=\{1\}\cup \{n\cdot 2n^3\} \cup \bigcup_{k=1}^{n-1}\left[(k-1)\cdot 2n^3+1, (k-1)\cdot 2n^3+k\cdot n\right],\\
Q(M_n)&:=\{1\}\cup \{n\cdot 2n^3\}\cup (M_n\setminus P(M_n)).
\end{align*}
We denote by $\mathcal{C}_{PQ}$ the class of structures $\{([1,n\cdot 2n^3],<,P,Q):n<\omega\}$.
\eej

Note that $P(M_n),Q(M_n)$ are predicates whose union is $M_n$, and whose intersection is given precisely by the maximal and minimal elements. We can describe the structure $M_n$ by considering it as the interval $[1,2n^4]$ divided into $n$ pieces of size $2n^3$, and such that in the $k$-th piece the predicate $P$ takes the first $k\cdot n$ elements (i.e., at least $n$ elements) while $Q$ takes the remaining $2n^3-k\cdot n$ elements. These structures provide an example of a class $\mathcal{C}_{PQ}$ whose infinite ultraproducts are all quasi-o-minimal (see Appendix \ref{Cpq-quasi-o-minimal}), and we will show in this section that $\mathcal{C}_{PQ}$ is not a weak o-asymptotic class. \\

It is easy to show that $\mathcal{C}_{PQ}$ is not an o-asymptotic class, as any decomposition of $M_n$ into a fixed number $\ell$ of intervals will necessarily contain arbitrarily large segments of elements in $P$ and $Q$, contradicting Lemma \ref{lemmaalternateintervals}. We show now that $\mathcal{C}_{PQ}$ does not even satisfy the weak version given in Definition \ref{def-weak-o-class}.

\bp \label{Cpq-not-o-class} The class $\mathcal{C}_{PQ}$ is not a weak o-asymptotic class.
\ep

\bdem If the class $\mathcal{C}_{PQ}$ were weak o-asymptotic, then for the formula $P(x)$ there would exist a constant $C>0$, a natural number $\ell$, tuples $(\mu_1,\ldots,\mu_\ell)\in [0,1]^\ell$ and tuples $\ov{c}^n=(c_0^n,\ldots,c_\ell^n)\in M_n^\ell$ such that for all $n$ and all $i\leq \ell$, either $\mu_i=0$ and $|P(M_n)\cap [c_{i-1}^n,c_i^n]|\leq C$ or $\mu_i>0$ and \[\left||P(M_n)\cap [c_{i-1}^n,c_i^n]|-\mu_i|[c_{i-1}^n,c_i^n]|\right|\leq C|[c_{i-1}^n,c_i^n]|^{1/2}.\]

If we analyze the first $\ell+1$ $Q$-intervals in $M_n$, by the pigeonhole principle there is one interval $I_i=[c_{i-1},c_i]$ in the decomposition that contains points of at least two different $Q$-intervals. So $I_i$ contains at least one full $P$-interval, and we have $|P(M_n)\cap I_i|\geq n$. Thus, by taking $n>C$ we conclude that $\mu_i>0$. \\

Suppose now that $n>\left(\frac{2C}{\mu_i}\right)^2$. Thus, we have both $\frac{\mu_i}{2}|I_i|>C|I_i|^{1/2}$ and $n>\frac{2}{\mu_i}$. If we show $\frac{|P(M_n)\cap I_i|}{|I_i|}<\frac{\mu_i}{2}$, then we would have \[|P(M_n)\cap [c_{i-1},c_i]|<\frac{\mu_i}{2}|[c_{i-1},c_i]|<\mu_i\cdot |[c_{i-1},c_i]|-C|[c_{i-1},c_i]|^{1/2},\] which would contradict condition (1) of the definition of weak o-asymptotic classes.\\

The maximal possible value for the fraction $\frac{|P(M_n)\cap I_i|}{|I_i|}$ occurs when the interval $I_i$ starts with a full $P$-interval and ends with another full $P$-interval. Thus, we may assume that $I_i$ is of the form \[I_i=P_r\cup Q_r\cup P_{r+1}\cup Q_{r+1}\cup \cdots \cup P_{r+t}\cup Q_{r+t}\cup P_{r+t+1}\]
where $P_s,Q_s$ denote the $P$-interval and the $Q$-interval in the $s$-th piece of $M_n$, respectively. In this case, since each $P$-interval in $M_n$ contains at most $n^2$ elements, we have
\begin{align*}
\dfrac{|P(M_n)\cap I_i|}{|I_i|}&\leq \dfrac{|P_r|+\cdots + |P_{r+t+1}|}{|P_{r+1}|+|Q_{r+1}|+\cdots + |P_{r+t}|+|Q_{r+t}|}=\dfrac{n^2(t+1)}{2n^3\cdot t}\leq \frac{1}{n}<\dfrac{\mu_i}{2}.
\end{align*}Therefore, for sufficiently large $n$, we obtain a contradiction.
\edem

\subsection{Classes of ordered graphs}

The idea behind the definition of o-asymptotic classes is to meld ideas of 1-dimensional asymptotic classes with some versions of o-minimality. So, it is expected to have examples of ultraproducts of o-asymptotic classes that are neither NIP nor simple. However, the assumptions of uniform distribution and definability on o-asymptotic classes impose some structural restrictions in their ultraproducts. For instance, in the following result, we show that the random graph cannot appear as the reduct of an ultraproduct of o-asymptotic classes. This contrasts with the example of the class of Paley graphs, which form a 1-dimensional asymptotic class and whose ultraproducts are elementarily equivalent to the random graph.

\bp There is no ultraproduct $M$ of an o-asymptotic class in the language $L=\{R,<\}$ such that the reduct $(M,R)$ satisfies the theory of the random graph.
\ep
\bdem Let $M$ be an infinite ultraproduct of structures in an o-asymptotic class $\mathcal{C}$ with language $L=\{<,R\}$, where $R$ is a binary irreflexive symmetric relation, and suppose that $(M,R)$ is a model of the theory of the random graph. Consider the $\omega$-type
\[p((x_s)_{s<\omega};y):=\{S^m(x_{2s})Ry \wedge \neg S^m(x_{2s+1})R y : s,m<\omega\}\cup\{x_{2s}<x_{2s+1}<x_{2s+2}:s<\omega\},\]where $S^m(z)$ represents the $m$-th succesor of $z$.\\

Using the properties of the random graph, it is easy to check that $p((x_s)_{s<\omega};y)$ is finitely satisfiable. Namely, given a finite subset $\Gamma_0$ of $p((x_s)_{i<\omega},y)$, let $m'=\max\{m:S^m(x_{2s}Ry\wedge \neg S^m(x_{2s+1})Ry\in \Gamma_0 \text{ for some $s$}\}$ and $k_0=\max\{s:x_s \text{ is mentioned in some formula in }\Gamma_0\}$. Then we can find elements $a_0,\ldots,a_k$ in $M$ such that $S^{m'}(a_s)<a_{s+1}$, and given the finite sets $X=\{S^m(a_{2s}):0\leq m\leq m', i\leq k\}$ and $Y=\{S^m(a_{2s+1}):0\leq m\leq m', s\leq k\}$, by the properties of the random graph there is an element $b$ that is adjacent with all elements in $X$ and non-adjacent with all elements in $Y$. Hence, $(a_s)_{s\leq k},b\models \Gamma_0$. \\

By compactness and $\aleph_1$-saturation, there are elements $(a_s:s<\omega), b$ in $M$ realizing $p$. Thus, we can take $\varphi(x,b):=xRb$ and define the infinite convex sets $I_s=\{S^m(a_{2s}):i<\omega\}$ and $J_s=\{S^m(a_{2s+1}):m<\omega\}$. By Lemma \ref{lemmaalternateintervals}, since $I_s\subseteq \varphi(M,b)$ and $J_s=M\setminus \varphi(M,b)$, we conclude that $\mathcal{C}$ is not an o-asymptotic class. 
\edem

\bp The class $\mathcal{OG}_k$ of linearly ordered graphs with fixed degree $k\geq 1$ is not an o-asymptotic class.
\ep
\bdem For every $n<\omega$, define the ordered graph $A_n$ with vertices $[1,\ldots,2n]=\mathbb{Z}/(2n)\mathbb{Z}$ and every vertex $i$ is adjacent to the vertices $i+n, i+1, i+2,\ldots, i+(k-1)$. Consider the graph $G_n=\underbrace{A_n\oplus \cdots \oplus A_n}_{\text{$n^2$ times}}$ formed by putting consecutively $n^2$ copies of the graph $A_n$. Let $\varphi(x)$ be the formula $\varphi(x):=\forall y\left(xRy\rightarrow x<y\right)$. Note that the formula $\varphi(x)$ defines in each copy $A_n$ of $G_n$ the set $[1,n]$. Hence, in an infinite ultraproduct $G$ of the graphs $G_n$ there will be infinite convex sets $\langle I_s<J_s<I_{s+1}:s<\omega\rangle$ with $I_s\subseteq \varphi(G)$ and $J_s\subseteq G\setminus \varphi(G)$ for every $s<\omega$, and we can conclude that $\mathcal{OG}_k$ is not an o-asymptotic class by Lemma \ref{lemmaalternateintervals}. 
\edem 

\brmk The previous proof shows in particular that several classes of linearly ordered graphs that are Ramsey classes  but  not o-asymptotic classes. In Section \ref{vs-lex-order} we give more examples of this phenomenon using finite vector spaces.\ermk

\subsection{Linear orders on classes of finite trees.}
\bd A \emph{finite tree} is a finite partial order $(T,\preceq)$ that has a unique $\preceq$-minimal element $\rho$, and such that for every $a\in T$, the set $T_{\preceq a}:=\{x\in T:x\preceq a\}$ is linearly ordered by $\preceq$. We say that a finite tree $T$ is a \emph{$k$-tree} if every element $a$ of $T$ has at most $k$ $\preceq$-successors, i.e., the set $S_{\preceq}(a)=\{x\in T: T\models \forall y(a\prec y\preceq x \rightarrow x=y)\}$ has cardinality less than or equal to $k$.

\ed

It is well-known that every finite $k$-tree $T$ can be seen as a subset of the set of finite sequences of elements in $\{1,\ldots,k\}$ of length at most $n$ (for some $n$), with $f\preceq g$ if and only if $f\subseteq g$. With this identification, there are two  natural ways to define linear orders on $T$:

\bd [Lexicographical order on a tree] Let $T\subseteq [k]^{\leq n}$ be a finite $k$-tree. We write $f<_{lex}g$ if and only if $f\subsetneq g$ or for some $i\leq \operatorname{length}(f)$, $f\upharpoonright_i=g\upharpoonright_i$ and $f(i)<g(i)$.
\ed
\bd [Level-based order on a tree] Let $T\subseteq [k]^{\leq n}$ be a finite $k$-tree. We define the \emph{level-based order $<_{lev}$} inductively on the length of the sequences, as follows:
\bitem
\item The empty sequence $\Lambda$ is the minimum of $<_{lev}$.
\item For sequences of length at most $1$ we impose
\[(1) <_{lev} (2) <_{lev} \cdots <_{lev} (k).\]
\item For arbitrary $f=(a_1,\ldots,a_{s+1}),g=(b_1,\ldots,b_{t+1})\in T$, $f<_{lev} g$ if and only one of the following three cases hold: either $s<t$, or $s=t$ and $f\upharpoonright_s=(a_1,\ldots,a_s)<_{lev}(b_1,\ldots,b_t)=g\upharpoonright_t$, or $f\upharpoonright_s=g\upharpoonright_t$ and $(a_{s+1})<_{lev}(b_{s+1})$.\eitem
\ed

\bej \label{example-lex-lev-order} Consider the tree $T=\{\Lambda,1,2,3,11,12,21,31,32,111\}\subseteq [3]^{\leq 3}$. The lexicographical order on $T$ is given by\[\Lambda<_{lex}1<_{lex}11<_{lex}111<_{lex}12<_{lex}2<_{lex}21<_{lex}3<_{lex}31<_{lex}32,\]
while the level-based order is \[\Lambda <_{lev} 1 <_{lev} 2 <_{lev} 3 <_{lev} 11 <_{lev} 12 <_{lev} 21 <_{lev} 31 <_{lev} 32 <_{lev} 111.\]
\eej
\bd Let us denote by $\mathcal{T}_{k,lex}$ the class of finite $k$-trees, equipped with the lexicographical order, and by $\mathcal{T}_{k,lev}$ the class of finite $k$-trees, equipped with the level-based order.
\ed
\bp The classes $\mathcal{T}_{k,lev}$ and $\mathcal{T}_{k,lex}$ are not o-asymptotic classes.\ep
\bdem
Consider the formula $\varphi(x,y)=x\preceq y$ (where $\preceq$ is the partial order on the tree) and let $T_{k,n}$ be the complete $k$-tree of height $n$. Choose $b_n$ to be the maximal element of the linear order, which in both cases is $b_n=(k,k,\ldots,k)\in [k]^n$.  Then, $\varphi(T_{k,n},b_n)$ consists of $n+1$ elements $a_{n,0}=\Lambda<a_{n,1}<\ldots<a_{n,n}=b_n$ and the $<$-interval $[a_i,a_{i+1}]$  contains $k^{i+1}+1$ elements (respectively $(k-1)\cdot k^{n-i-1}$ elements) of the structure $(T_{k,n},<_{lev})$ (resp. $(T_{k,n},<_{lex})$) which are not in $\varphi(T_{k,n},b_n)$. By taking an ultraproduct $T$ of the trees $T_{k,n}$ with respect to a non-principal ultrafilter $\mathcal{U}$ on the set of indices $n$, the set $X=\varphi(M,\max M)$ is an infinite set with measure zero in every infinite subinterval. Thus, by Lemma \ref{lemmameasurezero}, we conclude that neither $\mathcal{T}_{k,lev}$ nor $\mathcal{T}_{k,lex}$ are o-asymptotic classes.
\edem

\subsection{Vector spaces with the lexicographical order} \label{vs-lex-order} Let $p$ be a fixed prime and consider the finite vector spaces $\{\mathbb{F}_p^n:n\geq 1\}$ over $\mathbb{F}_p$. By identifying $\mathbb{F}_p$ with $\mathbb{Z}/p\mathbb{Z}$, we can use the order on the cyclic group defined in Section \ref{examples-O-classes} to induce a lexicographical order on $\mathbb{F}_p^n$. There are three classes to consider depending on whether we fix the dimension or the base field.

\bd Let $\mathbb{P}$ be the set of prime numbers. We define:\benum
\item The class of ordered finite vector spaces $\C_{o-vs}:=\{(\mathbb{F}^n_p,+,<_{lex}):p\in\mathbb{P},n<\omega\}.$
\item The class of ordered finite vector spaces over $\mathbb{F}_p$, $\C_{o-vs,\,p}:=\{(\mathbb{F}^n_p,+,<_{lex}):n<\omega\}.$
\item The class of ordered finite vector spaces of dimension $k$ over fields of prime cardinality, $\C_{o-vs,\,\dim=k}=\{(\mathbb{F}^k_p,+,<_{lex}):p\in\mathbb{P}\}.$
\eenum
\ed
To start, note that for every odd prime $p$ we can put $q=\frac{p-1}{2}$, so that $q$ is the maximal element of $(\mathbb{F}_p,+,\ci{<})$.
\bp The class $\mathcal{C}_{o-vs,\,p}$ is not an o-asymptotic class.
\ep
\bdem Suppose $p\geq 5$. Then, in $(\mathbb{F}_p,+,<)$ the inequation $-1<2\cdot t$ has exactly $\frac{p+1}{2}$ solutions (a unique solution for each equation $2\cdot t=a$ with $a=0,1,\ldots,q$), and $t=-1$ is not a solution for the inequation. Also, the equation $2t=-1$ has exactly one solution, namely $t=q$.\\

For a fixed $n<\omega$, let $b_n=(-1,-1,\ldots,-1)\in\mathbb{F}_p^n$ and consider the formula $\varphi(x;b_n):=x+x>b_n$. For $s=1,\ldots,n-1$, we can define the sets 
\begin{align*}
Y_{s,n}&=\left\{(x_1,\ldots,x_n)\in\mathbb{F}_p^n: x_1=x_2=\cdots=x_{s-1}=q \text{ and }-1<2\cdot x_s\right\}\\
&= \bigcup_{\substack{a\in\mathbb{F}_p\\ 2\cdot a>-1}} \left(\{(\underbrace{q,\ldots,q}_{s-1 \text{ times}},a)\}\times \mathbb{F}_p^{n-s}\right).
\end{align*}
By definition of the lexicographical order, we have $Y_{1,n}<_{lex}Y_{2,n}<_{lex}\cdots<_{lex}Y_{n,n}$ and $\varphi(\mathbb{F}_p^n;b_n)=Y_{1,n}\cup \cdots \cup Y_{n,n}$. For $s\leq \frac{n}{2}$ we can define the intervals \[I_{s,n}=\left\{(\underbrace{q,\ldots,q}_{s-1 \text{ times}},0)\right\}\times \mathbb{F}_{p}^{n-s}, \hspace{0.5cm}J_{s,n}=\left\{(\underbrace{q,\ldots,q}_{s-1 \text{ times}},-1)\right\}\times \mathbb{F}_{p}^{n-s}.\]

In the lexicographical order, we have $I_{1,n}<_{lex}J_{1,n}<_{lex}\cdots <_{lex} I_{\frac{n}{2},n}<J_{\frac{n}{2},n}$, and we have $I_{s,n}\subseteq Y_{s,n}\subseteq \varphi(\mathbb{F}_p^n;\ov{b}_n)$, $J_{s,n}\subseteq \mathbb{F}_p^n\setminus \varphi(\mathbb{F}_p^n;\ov{b}_n)$. In any infinite ultraproduct $F$ of the structures $(\mathbb{F}_p^n,+,<_{lex})$ there will be infinite intervals $\langle I_s,J_s:s<\omega\rangle$ satisfying the hypothesis of Lemma \ref{lemmaalternateintervals}. Hence, $\mathcal{C}_{o-vs,\,p}$ is not an o-asymptotic class. \\

The previous argument also works when $p=3$, but using the formula $\varphi(x):=x+x>0$. Finally, when $p=2$, we can define $M_n=\mathbb{F}_2^{n^2}$ and $b_n\in \mathbb{F}_2^{n^2}$ to be the element that has $1$ only in the $k\cdot n$-th coordinates, for $k=1,2,\ldots,n$. If we take $\varphi(x,b):=x+b<x$, then $M_n\models \varphi(x,\ov{b})$ for $x=(x_1,\ldots,x_{n^2})$ if and only if $x_{k\cdot n}=1$ for each $k=1,2\ldots,n$. Thus, if $M=\prod_\mathcal{U} M_n$ is an infinite ultraproduct of the structures $M_n$, and $b=[b_n]_\mathcal{U}$, we can conclude using Lemma \ref{lemmameasurezero} that $\mathcal{C}_{o-vs,\,2}$ is not an o-asymptotic class, as for every infinite interval $[\alpha,\beta]$ of $M$ we have \[\operatorname{meas}_{[\alpha,\beta]}(\varphi(M,\ov{b}))\leq \lim_{n\to\infty} \dfrac{|\varphi(M_n,\ov{b}_n)|}{|M_n|}=\lim_{n\to\infty}\dfrac{|\mathbb{F}_2^{n^2-n}|}{\mathbb{F}_2^{n^2}}=\lim_{n\to\infty}\frac{1}{2^n}=0.\qedhere\]
\edem

The previous result also shows that $\mathcal{C}_{o-vs}$ is not an o-asymptotic class. We now focus on the class $\mathcal{C}_{o-vs,\dim=k}$. Since the dimension in the class is fixed to be $k$, we can define for $i=1,\ldots,k$ the elements $e_i:=(0,\ldots,0,1,0,\ldots,0)$ and $q_i:=(0,\ldots,0,q,0,\ldots,0)$.

\bl The elements $e_1,\ldots,e_k,q_1,\ldots,q_k$ are uniformly definable over the empty set in every structure of the class $\mathcal{C}_{o-vs,\dim=k}$, and so are the projection functions $\pi_1,\ldots,\pi_k$.
\el
\bdem We will define recursively the elements $e_k,e_{k-1},\ldots,e_1$ and $q_k,q_{k-1},\ldots,q_1$. To start, note that $e_k=(0,\ldots,0,1)$ is precisely the successor of $0=(0,\ldots,0)$ in the lexicographical order. Also, note that $q_k=(0,\ldots,0,q)$ is precisely the minimum element $x$ satisfying $x>0$ and $x+e_k<x$.\\

Suppose now that $e_{k},\ldots,e_{i+1}$ and $q_{k},\ldots,q_{i+1}$ have been already defined. Notice that $q_{i+1}+\ldots+q_k=(0,\ldots,0_i,q,q,\ldots,q)$, and if $(0,\ldots,0_i,q,q,\ldots,q)<(a_1,\ldots,a_i,\ldots,a_k)$ we have $a_j>0$ for some $j\leq i$, as it is impossible to have $a_j>q$ for some $j=i+1,\ldots,k$. Hence, $(0,\ldots,0,1,-q,-q,\ldots,-q)\leq (a_1,\ldots,a_i,\ldots,a_k)$. Therefore, if $S$ denotes the successor function in the lexicographical order, we have \begin{align*}
S(q_{i+1}+\cdots+q_k)+(q_{i+1}+\ldots+q_k)&=(0,\ldots,0,1_i,-q,-q,\ldots,-q)+(0,\ldots,0,0_i,q,q,\ldots,q)\\ &=(0,\ldots,0,1_i,0,\ldots,0)=e_i.
\end{align*}
Consider now the minimal element $x$ satisfying $x>0$ and $x+e_i<x$. Note that if $a=(a_1,\ldots,a_k)$ satisfies $a>0$ and $a+e_i<a$ we have either $a_j\geq 0$ for some $j< i$, or $a=(0,\ldots,0,a_i,a_{i+1},\ldots,a_k)$ and $a+e_i=(0,\ldots,0,a_i+1,a_{i+1},\ldots,a_k)<(0,\ldots,0,a_i,a_{i+1},\ldots,a_k)$ implies $a_i=q$. In either case, $(0,\ldots,0,q,-q,\ldots,-q)\leq_{lex}(a_1\ldots,a_k)=a$. Hence, $q_i$ can be defined from $e_i,q_{i+1},\ldots,q_k$ as the unique element $z$ satisfying
\[\exists y\left(y>0\wedge y+e_i<y\wedge \forall x(x>0\wedge x+e_i<x\rightarrow y\leq x)\wedge z=y+q_{i+1}+\ldots+q_k\right).\]

Finally, notice that for $i=1,2,\ldots,k$, the $x_i$-th axis can be defined by the formula $A_i(x)=\bigwedge_{j\neq i}(x-q_j<x<q_j)$. Hence, for every $x\in \mathbb{F}_p^k$, the tuple of projections $(\pi_1(x),\ldots,\pi_k(x))$ is the unique tuple $(y_1,\ldots,y_k)$ satisfying $(\mathbb{F}_p^k,+,<_{lex})\models x=y_1+\cdots + y_k+\bigwedge_{i=1}^k A_i(y_i)$.
\edem

\bp The class $\mathcal{C}_{o-vs,\dim=k}$ is not an o-asymptotic class.
\ep
\bdem Consider the formula $\varphi(x):=\pi_{k}(x)=0$, which defines in each $\mathbb{F}_p^k$ the set \[X_p=\{(x_1,\ldots,x_k)\in\mathbb{F}_p^k:x_k=0\}.\] In every infinite ultraproduct $F$ of the structures $\mathbb{F}_p^k$ (with respect to a non-principal ultrafilter on $\mathbb{P}$) $\varphi(F)$ is infinite and has measure zero in every infinite interval of $F$. Hence, by Lemma \ref{lemmameasurezero}, $\mathcal{C}_{o-vs,\dim=k}$ is not an o-asymptotic class.
\edem

\section{Appendix}\label{Appendix}
\subsection{Uniform quantifier elimination for the class of finite linear orders.} \label{uniformqe} \ \\

We start this appendix by recalling a version of quantifier elimination for classes of structures.
\bd Let $\mathcal{C}$ be a class of $L$-structures. We say that the class has \emph{uniform quantifier elimination} if for every $L$-formula $\theta(\ov{y})$ there is a quantifier-free formula $\psi(\ov{y})$ such that for every $M\in\mathcal{C}$, $M\models \forall \ov{y}(\theta(\ov{y})\leftrightarrow \psi(\ov{y}))$.\\

This is equivalent to say that the theory $\operatorname{Th}(\mathcal{C})$, which is the common theory of all the structures in the class $\mathcal{C}$, has quantifier elimination.
\ed

In this appendix we will show that both the class $\mathcal{C}_{ord}$ of Example \ref{example-Cord} and the class $\mathcal{C}_{\text{$k$-col}}$ of Example \ref{example-Ckcol} have uniform quantifier elimination in certain suitable languages. 
\bl \label{UQEcord} The class $\mathcal{C}_{ord}$ has quantifier elimination in the language $L'=\{<,\min,\max,S,S^{-1}\}$  
\el
\bdem Recall that every structure in $\mathcal{C}_{ord}$ has a canonical expansion to an $L'$-structure by interpreting the constant symbols $\min,\max$ to be the minimum and maximum elements, and the unary functions $S,S^{-1}$ to be respectively the successor and the predecessor functions, defining also $S(\max)=\max$ and $S^{-1}(\min)=\min$.\\

In this language, the atomic formulas have the form $\tau_1=\tau_2, \tau_1<\tau_2$ for terms $\tau_1,\tau_2$ in $L'$. Since $<$ is a linear order, we have that $\neg(\tau_1=\tau_2)\equiv \tau_1<\tau_2\vee \tau_2<\tau_1$ and $\neg(\tau_1<\tau_2)\equiv \tau_1=\tau_2 \vee \tau_1>\tau_2$. Therefore, every \emph{primitive existential formula}\footnote{That is, an existential formula quantifying a variable in a conjunction of atomic formulas and negations of atomic formulas. It is well-known that in order to show quantifier elimination it is enough to prove that every primitive existential formula is equivalent to a quantifier-free one.} can be assumed (possibly after distributing the existential quantifier in disjunctions) to have the form

\[\theta(\ov{y}):=\exists x\left(\bigwedge_{i=1}^k S^{m_i}x=\sigma_i(\ov{y}) \wedge \bigwedge_{j=1}^\ell \tau_{1,j}(\ov{y})<S^{n_j}x<\tau_{2,j}(\ov{y})\right)\]
for some integers $m_1,\ldots,m_k,n_1,\ldots,n_\ell$ and some $L'$-terms $\sigma_i,\tau_{1,j},\tau_{2,j}$ (if necessary, we could take $\tau_{1,j}=\min$ or $\tau_{2,j}=\max$). We now analyze each of the conjuncts separately. 

For $S^{m_i}x=\sigma_i(\ov{y})$, we have:
\bitem
\item If $m_i\geq 0$, $S^{m_i}x=\sigma_i(\ov{y})$ is equivalent to \[(S^{-m_i}\sigma_i(\ov{y})\leq x\leq S^{-m_i}\sigma_i(\ov{y})) \vee (\sigma_i(\ov{y})=\max \wedge\  S^{-m_i}(\max) \leq x).\]
\item If $m_i<0$, $S^{m_i}x=\sigma_i(\ov{y})$ is equivalent to
\[(S^{-m_i}\sigma_i(\ov{y}) \leq x\leq S^{-m_i}\sigma_i(\ov{y}) ) \vee (\sigma_i(\ov{y})=\min \wedge\  \min\leq x\leq S^{-m_i}(\sigma_i(\ov{y}))).\]
\eitem
Similarly, the formula $\tau_{1,j}(\ov{y})<S^{n_j}x$ is equivalent to \[(S^{-n_j+1}(\tau_{1,j}(\ov{y}))\leq x\leq \max) \vee  (\max>\tau_{1,j}(\ov{y})\wedge S^{-n_j}(\max)\leq x\leq \max),\]
and the formula $S^{n_j}x<\tau_{2,j}(\ov{y})$ is equivalent to \[(\min\leq x\leq S^{-n_j-1}(\tau_{2,j}(\ov{y}))) \vee  (\min<\tau_{2,j}(\ov{y})\wedge \min\leq x \leq S_{-n_j}(\min)).\]

Thus, after distributing further the existential quantifier on disjunctions, there are $L'$-terms $\langle\eta_{i,r,s}(\ov{y}): i=1,2; r\leq k', s\leq \ell'\rangle$ such that $\theta(\ov{y})$ is equivalent to 
\begin{align*}
&\exists x\left(\bigvee_{r\leq k'}\bigwedge_{s\leq \ell'} \eta_{1,r,s}(\ov{y})\leq x\leq \eta_{2,r,s}(\ov{y})\right)\equiv \bigvee_{r\leq k'}\exists x \left( \bigwedge_{s\leq \ell'} \eta_{1,r,s}(\ov{y})\leq x\leq \eta_{2,r,s}(\ov{y})\right)\\
&\equiv \bigvee_{r\leq k'}\bigwedge_{s,s'\leq \ell'}\eta_{1,r,s}(\ov{y})\leq \eta_{2,r,s'}(\ov{y}).
\end{align*}
This finishes the proof that $\mathcal{C}_{ord}$ has quantifier elimination in the language $L'$.
\edem
\bl \label{UQE-k-col} The class $\mathcal{C}_{k-col}=\{([1,n],<,P_1,\ldots,P_k):n<\omega\}$ has uniform quantifier elimination in the language $L'=\{<,S,S^{-1},P_1\}$.
\el
\bdem First, notice that the predicates $P_i$ can all be expressed in the language $L'$, as for structures $M\in\mathcal{C}_{k-col}$ we have  $M\models P_i\Leftrightarrow M\models S^{i-1}(\min)\leq x \wedge P_1(S^{i-1}(x)).$
Now, suppose $\theta(\ov{y})$ is a primitive existential formula in the language $L'$, which as in the previous proof can be assumed to have the form
\begin{align*}
&\theta_1(\ov{y}):=\exists x\left(\psi(\ov{y})\wedge P_1(x)\wedge \bigwedge_{j=1}^\ell \tau_{1,j}(\ov{y})\leq x\leq \tau_{2,j}(\ov{y}) \right), \text{\ \ or}\\
&\theta_2(\ov{y}):=\exists x\left(\psi(\ov{y})\wedge \neg P_1(x)\wedge \bigwedge_{j=1}^\ell \tau_{1,j}(\ov{y})\leq x\leq \tau_{2,j}(\ov{y}) \right)
\end{align*}
for some quantifier-free formula $\psi(\ov{y})$, and $L'$-terms $\tau_{1,j}(\ov{y}),\tau_{2,j}(\ov{y})$. In the first case, we simply have
\[\theta_1(\ov{y})\equiv \bigwedge_{j=1}^\ell \left(S^k(\tau_{1,j}(\ov{y}))\leq \tau_{2,j}(\ov{y}) \vee \bigvee_{1\leq r<s\leq k}(P_s(\tau_{1,j}(\ov{y}))\wedge P_r(\tau_{2,j}(\ov{y}))\right).\]
In the second case, since $\neg P_1(x)$ is equivalent to $\displaystyle{\bigvee_{2\leq t\leq k} P_1(S^{-t}x)}$, we have 
\begin{align*}
\theta_2(\ov{y})
&\equiv\exists x\left(\bigvee_{2\leq t\leq k} P_1(S^{-t}x)\wedge \bigwedge_{j=1}^\ell \tau_{1,j}(\ov{y})\leq x\leq \tau_{2,j}(\ov{y}) \right)\\
&\equiv\bigvee_{2\leq t\leq k}\exists x\left(P_1(S^{-t}x)\wedge \bigwedge_{j=1}^\ell S^{-t}(\tau_{1,j}(\ov{y}))\leq S^{-t}(x)\leq S^{-r}(\tau_{2,j}(\ov{y}) )\right)\\
&\equiv \bigvee_{2\leq t\leq k}\bigwedge_{j=1}^\ell \left(S^{k-t}(\tau_{1,j}(\ov{y}))\leq S^{-t}(\tau_{2,j}(\ov{y}) )\vee \bigvee_{1\leq r<s\leq k} P_1(S^{-t-(s-1)}(\tau^1_j(\ov{y})))\wedge P_1(S^{-t-(r-1)}(\tau_{2,j}(\ov{y})))\right).
\end{align*}This finishes the proof.
\edem

\subsection{Infinite ultraproducts of the class $\mathcal{C}_{PQ}$ are quasi-o-minimal.} \label{Cpq-quasi-o-minimal} 

Let $\mathcal{U}$ be a non-principal ultrafilter on $\Nat$ and consider the structure $M^*=\prod_{\mathcal{U}}M_n$ where the finite structures $M_n$ are members of the class $\mathcal{C}_{PQ}$, defined in Example \ref{defCp}. In the structure $M^*$ we can define the unary functions $P_+,P_-,Q_+,Q_-$ in the following way:

\begin{align*}
P_+(x)=y\text{ if and only if }[P(x)\wedge x=y] \vee [Q(x)\wedge P(y)\wedge \forall z(x<z<y\rightarrow Q(z))],\\
P_-(x)=y\text{ if and only if }[P(x)\wedge x=y] \vee [Q(x)\wedge P(y)\wedge \forall z(y<z<x\rightarrow Q(z))],\\
Q_+(x)=y\text{ if and only if }[Q(x)\wedge x=y] \vee [P(x)\wedge Q(y)\wedge \forall z(x<z<y\rightarrow P(z))],\\
Q_-(x)=y\text{ if and only if }[Q(x)\wedge x=y] \vee [P(x)\wedge Q(y)\wedge \forall z(y<z<x\rightarrow P(z))].
\end{align*}
We can see $P$ and $Q$ as some kind of parity predicates (except for the minimum and maximum), so $P_+(x)$ is the first element that is bigger than or equal to $x$ with ``parity'' $P$, whereas $Q_-(x)$ is the first element obtained moving from $x$ to the left that has parity $Q$.\\

The following identities follow easily from the definitions:
\begin{align*}
&P_+(P_-(x))=P_-(x) & & P_-(P_+(x))=P_+(x)\\
&Q_+(Q_-(x))=Q_-(x) & & Q_-(Q_+(x))=Q_+(x)\\
&P_+(Q_-(x))=\begin{cases}P_+(x)&\text{if $x\in Q$}\\ S(Q_-(x))&\text{if $x\in P$}\end{cases} & &P_-(Q_+(x))=\begin{cases}P_-(x)&\text{if $x\in Q$}\\ S^{-1}(Q_+(x))&\text{if $x\in P$}\end{cases}\\
&Q_+(P_-(x))=\begin{cases}Q_+(x)&\text{if $x\in P$}\\ S(P_-(x))&\text{if $x\in Q$}\end{cases} & &Q_-(P_+(x))=\begin{cases}Q_-(x)&\text{if $x\in P$}\\ S^{-1}(P_+(x))&\text{if $x\in Q$}\end{cases}
\end{align*}

\bp The theory $T=\operatorname{Th}(M^*)$ has quantifier elimination in the extended language $L'=\{<,P,Q\}\cup\{\min,\max,S,S^{-1},P_+,P_-,Q_+,Q_-\}$.
\ep
\bdem
We use here a well-known criterion for quantifier elimination, that can be found in Theorem 5.4 of \cite{HensonMT}. Let $M,N$ be models of $T$, with $N$ $\omega$-saturated. Let $A=\langle a_1,\ldots,a_n\rangle$ be a common finitely generated substructure of $M$ and $N$ in the language $L'$, and let $b\in M\setminus A$. We aim to find an element $c\in N$ such that the map $\langle A,b\rangle\to \langle A,c\rangle$ is a partial embedding.\\

Since $b\not\in A$, we may assume without loss of generality that $a_1<b<a_2$ and the interval $(a_1,a_2)$ contains no element $a_i$, $i=1,2,\ldots,n$. Furthermore, since we have the successor and predecessor functions in $L'$, we may assume that $S^n(a_1)<b<S^{-n}(a_2)$ for every $n<\omega$. We denote this by $a_1<<b<<a_2$. There are essentially two cases to consider:\\

\noindent \emph{Case 1:} $\sigma(a_1)<<b<<\tau(a_1)$ or $\sigma(a_2)<<b<<\tau(a_2)$ for some $L'$-terms $\sigma,\tau$.\\

From the identities above we can conclude that the elements in $A$ can be obtained from $\{a_1,\ldots,a_n\}$ by iteratively alternating the function $P_+$ with $Q_+$ (and the function $P_-$ with $Q_-$) on each $a_i$, and then closing under successors and predecessors. Suppose for instance that we have $\tau(a_1)<<b<<P_+\tau(a_1)$ for some term $\tau$ that is a term of the form $Q_+P_+\cdots Q_+$. In this case, we know that the parity of $b$ is $Q$. Furthermore, we have
\begin{align*}
P_+(b)=P_+(\tau(a_1)) & &P_-(b)=S^{-1}(\tau(a_1))& &Q_-(b)=b=Q_+(b)
\end{align*}
Note that the set of formulas $\Gamma(x)=\{Q(x)\}\cup\{\tau(a_1)<S^{-n}(x),S^n(x)<P_+(\tau(a_1)):n<\omega\}$ is finitely consistent in $N$, as we can realize every finite set by an element of the form $S^k(\tau(a_1))$ for sufficiently large $k$. Since $N$ is $\omega$-saturated, there is an element $c\in N$ realizing $\Gamma(x)$, and for such element $c$ we would have $\tau(a_1)<<c<<P_+(\tau(a_1))$ and $P_+(c)=P_+(\tau(a_1)), P_-(c)=S^{-1}(\tau(a_1)), Q_-(c)=c=Q_+(c)$, thus obtaining a partial embedding $\langle A,b\rangle\to\langle A,c\rangle$.

The other possibilites can be analyzed in a similar way. \\

\noindent \emph{Case 2:} $\sigma(a_1)<<b<<\tau(a_2)$ for all $L'$-terms $\sigma,\tau$.\\

In this case we must also have $\sigma(a_1)<<\nu(b)<<\tau(a_2)$ for all $L'$-terms $\sigma,\tau,\nu$. If $b\in P$, the set $\Gamma_P(x)=\{P(x)\}\cup \{\sigma(a_1)<S^{-n}(\nu(x)),S^n(\nu(x))<\tau(a_2):\text{$\sigma,\tau,\nu$ $L'$-terms}, n<\omega\}$ is finitely consistent in $N$, as a finite subset $\Gamma_0$ can be witnessed by $P_+(Q_+(\sigma(a_1)))$, where $\sigma(a_1)$ is the largest term depending on $a_1$ mentioned in $\Gamma_0$. Hence, by $\omega$-saturation, there is a realization $c\in N$ of the set $\Gamma_P(x)$, and the map $\langle A,b\rangle\to\langle A,c\rangle$ is a partial embedding as all the inequalities hold and there are no non-trivial equalities to preserve. A similar argument can be used if $b\in Q$.
\edem
We now turn our attention to show that every ultraproduct of structures in the class $\mathcal{C}_{PQ}$ is \emph{quasi-o-minimal}, a weakening of the notion of o-minimality that was introduced by Belegradek, Peterzil and Wagner in \cite{BPW}.

\bd \label{def-quasi-o-min} A structure $M$ is \emph{quasi-o-minimal} if for every structure $N$ elementarily equivalent to $M$ the definable sets in one variable are finite boolean combinations of intervals with end-points in $N$ and $\emptyset$-definable sets.
\ed  
The main examples of quasi-o-minimal structures include the structure $(\mathbb{R},<,\mathbb{Q})$ (where $\mathbb{Q}$ is the realization of a unary predicate $Q(x)$) and the structure $(\mathbb{Z},+,-,0,1,<)$ of Presburger arithmetic, whose definable sets in one-variable are boolean combinations of intervals and the $\emptyset$-definable sets $D_m(x+t(\ov{y}))$ described in Section \ref{examples-O-classes}. \\

To show that every infinite ultraproduct of elements in the class $\mathcal{C}_{PQ}$ is quasi-o-minimal, we recall the following criterion:
\bh [Belegradek, Peterzil, Wagner - Theorem 1 in \cite{BPW}]\label{BPW-criterion}  A structure $M$ is quasi-o-minimal if and only if for every formula $\varphi(x,\ov{y})$ there is a formula $\chi(x,\ov{y},\ov{z})$ of the form \[\bigvee_i \varphi_i(x)\wedge \psi_i(\ov{y})\wedge \rho_i(x,\ov{z})\]such that the formulas $\rho_i(x,\ov{z})$ are conjunctions of formulas of the form $x=z,x<z,x>z$ and we have \[M\models \forall \ov{y}\,\exists \ov{z} \,\forall x\left(\varphi(x,\ov{y})\leftrightarrow \chi(x,\ov{y},\ov{z})\right).\]
\eh
\bl In ultraproducts of structures in the class $\mathcal{C}_{PQ}$, every formula $\varphi(x,y)$ in the language $L'$ is  equivalent to a positive combination of formulas of the form $P(x),Q(x),P(y),Q(y)$ and $x<\tau(y),x=\tau(y),x>\tau(y)$ where $\tau$ ranges over $L'$-terms.
\el
\bdem Since $L'$ contains only unary functions and the relation symbols $<,=$, every atomic formula in $L'$ using more than one variable has the form $\tau_1(x)=\tau_2(y)$, $\tau_1(x)<\tau_2(y)$ or $\tau_1(x)>\tau_2(y)$ for $L'$-terms $\tau_1,\tau_2$. So, by induction on terms, it is enough to prove the statement in the case where $\tau_1$ is a unary function in the language $L'$. It is clear when $\tau_1$ is $S$ or $S^{-1}$.

Suppose now $\tau_1=P_+$. Then we have:
\begin{align*}
P_+(x)=\tau_2(y)&\equiv (P(x)\wedge x=\tau_2(y))\vee (Q(x)\wedge P_-(Q_-(\tau_2(y))<x<y).\\
P_+(x)<\tau_2(y) &\equiv (P(x)\wedge x<\tau_2(y))\vee (Q(x)\wedge x<
P_-(S^{-1}(\tau_2(y)))).\\
P_+(x)>\tau_2(y) & \equiv (y<x) \vee (P_-(x)<x<y).
\end{align*}
By symmetry of the definitions of $P_+,P_-,Q_+,Q_-$, we can solve the other cases in a similar way.
\edem
\bp Every infinite ultraproduct $M$ of structures in the class $\mathcal{C}_{PQ}$ is quasi-o-minimal.
\ep
\bdem Consider a formula $\theta(x,\ov{y})$. By quantifier elimination in $L'$, $\theta(x,\ov{y})$ is equivalent to a formula of the form $\bigvee_{i}\eta_{i}(x,\ov{y})$, where each $\eta_{i}(x,\ov{y})$ is a conjunction of formulas that are either atomic formulas or negations of atomic formulas. Since we have a linear order and the only relation symbols are $=$ and $<$, the negation of an atomic formula can be written as a positive combination of atomic formulas, and by distributing the existential quantifiers in the disjunctions we may assume that $\eta_i(x,\ov{y})$ is already a positive conjunction of atomic formulas.\\

Each conjunction $\eta_i$ can be written as $\varphi_i(x)\wedge \psi_i(\ov{y})\wedge \rho_i'(x,\ov{y})$, where $\varphi_i(x)$ is the conjunction of all atomic formulas in $\eta_i(x,\ov{y})$ mentioning only the variable $x$, $\psi_i(\ov{y})$ is the conjunction of those mentioning only variables from $\ov{y}$, and $\rho_i'(x,\ov{y})$ is the conjunction of all atomic formulas mentioning both $x$ and at least one variable from $\ov{y}$. Since the language only contains unary functions, $\rho_i(x,\ov{y})$ is written as a conjunction of atomic formulas where only one variable from $\ov{y}$ appears at a time.\\

By the previous lemma, $\rho_i'(x,\ov{y})$ is equivalent to a positive combination of formulas of the form $P(x),Q(x)$ and $x=\tau(y),x>\tau(y),x<\tau(y)$ for an $L'$-term $\tau$ and a single variable $y$ from $\ov{y}$.  Let $\ov{z}=(z_1,\ldots,z_k)$ where $k$ is the number of different terms $\tau(y)$ appearing in $\rho'(x,\ov{y})$, and take $\rho_i(x,\ov{z})$ to be the formula obtained by replacing each term $\tau(y)$ by the $z_i$ corresponding to the enumeration. To finish the proof, let $\chi(x,\ov{y},\ov{z})=\bigvee_i\left(\varphi_i(x)\wedge \psi_i(\ov{y})\wedge \rho_i(x,\ov{z})\right)$. Thus, we have that each $\rho_i(x,\ov{z})$ is a conjunction of formulas of the form $x=z,x<z,x>z$, and we have that for all $\ov{a}\in M^{|\ov{y}|}$, the tuple $\ov{b}=(\tau_1(\ov{a}),\ldots,\tau_k(\ov{a}))$ realizes $M\models \forall x\left(\theta(x,\ov{a})\leftrightarrow \chi(x,\ov{a},\ov{b})\right)$.
 Hence, by Fact \ref{BPW-criterion}, the structure $M$ is quasi-o-minimal.
\edem

\begin{flushleft}
\begin{table}[H]
\begin{tabular}{ |>{\centering\arraybackslash}m{1.5cm} >{\arraybackslash}m{8cm} |}
\hline
\vspace{0.1cm}\includegraphics[width=11mm,height=10mm]{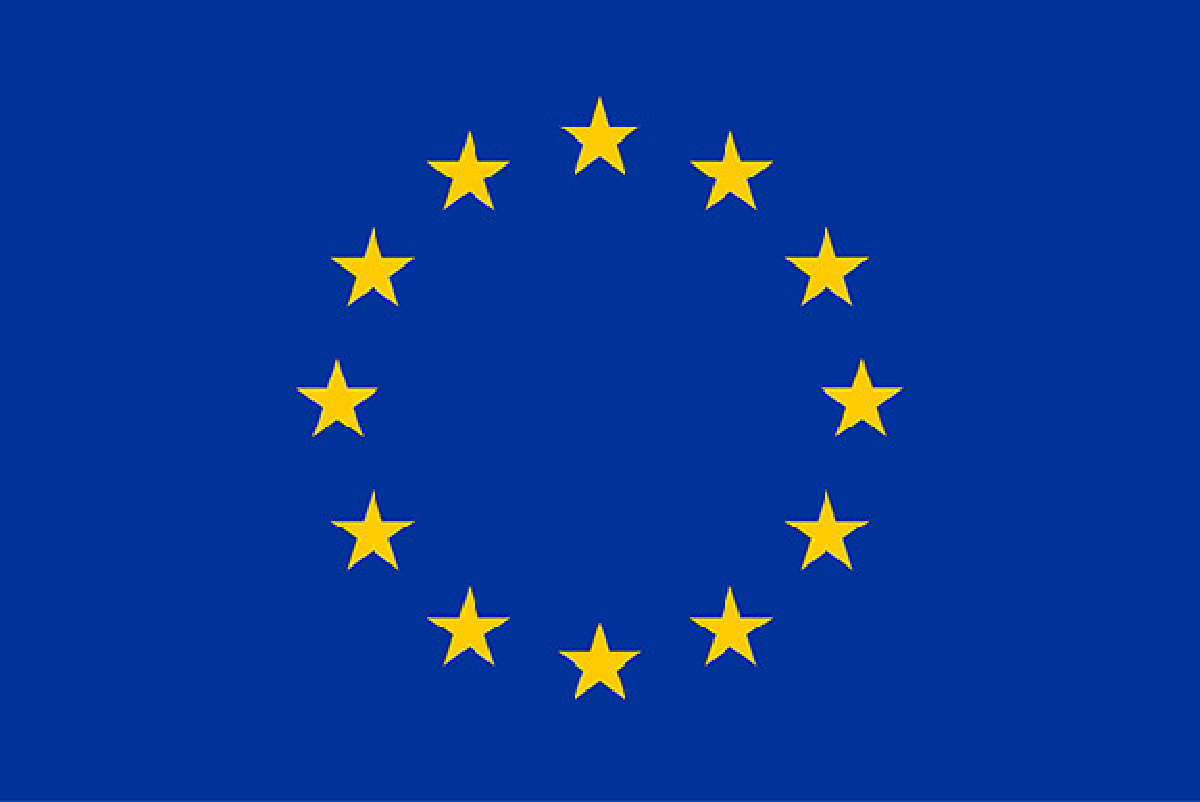} & \tiny{This project has received funding from the European Union's Horizon 2020 research and innovation programme under grant agreement No 656422} \\
\hline
\end{tabular}
\end{table}
\end{flushleft}

\end{document}